\documentclass[12pt, a4paper]{amsart}
\usepackage{latexsym}
\usepackage[mathcal]{eucal}
\usepackage{amssymb}
\usepackage{indentfirst}
\usepackage{amsmath,amsthm}
\usepackage{color,amsfonts}

\usepackage{amscd}

\allowdisplaybreaks[4] 

\begin{document}

\newcommand{\Norma}{\operatorname{N}}
\newcommand{\Ker}{\operatorname{Ker}}
\newcommand{\Ima}{\operatorname{Im}}
\newcommand{\con}{\operatorname{C}}
\newcommand{\gr}{\operatorname{gr}}
\newcommand{\Lie}{\operatorname{Lie}}
\newcommand{\Ad}{\operatorname{Ad}}
\newcommand{\dif}{\operatorname{d}\!}
\newcommand{\U}{\operatorname{U}}
\newcommand{\Z}{\operatorname{Z}}
\newcommand{\pro}{\operatorname{pro}}
\newcommand{\Hom}{\operatorname{Hom}}
\newcommand{\End}{\operatorname{End}}
\newcommand{\Tr}{\operatorname{Tr}}
\newcommand{\Ev}{\operatorname{Ev}}
\newcommand{\Span}{\operatorname{Span}}
\newcommand{\E}{\operatorname{E}}
\newcommand{\degree}{\operatorname{deg}}
\newcommand{\Invf}{\operatorname{Invf}}
\newcommand{\Inv}{\operatorname{Inv}}
\newcommand{\oc}{\operatorname{c}}
\newcommand{\Socle}{\operatorname{Socle}}
\newcommand{\oL}{\operatorname{L}}
\newcommand{\GL}{\operatorname{GL}}
\newcommand{\SL}{\operatorname{SL}}
\newcommand{\oH}{\operatorname{H}}

\newcommand{\g}{\mathfrak g}
\renewcommand{\k}{\mathfrak k}
\newcommand{\h}{\mathfrak h}
\newcommand{\p}{\mathfrak p}
\newcommand{\q}{\mathfrak q}
\renewcommand{\a}{\mathfrak a}
\renewcommand{\b}{\mathfrak b}
\renewcommand{\c}{\mathfrak c}
\newcommand{\n}{\mathfrak n}
\renewcommand{\u}{\mathfrak u}
\newcommand{\e}{\mathfrak e}
\newcommand{\f}{\mathfrak f}
\renewcommand{\l}{\mathfrak l}
\renewcommand{\t}{\mathfrak t}

\newcommand{\iZ}{\mathbb{Z}}
\newcommand{\C}{\mathbb{C}}
\newcommand{\R}{\mathbb R}
\renewcommand{\S}{\mathbf S}
\newcommand{\M}{\mathbf{M}}
\newcommand{\V}{\mathbf{V}}
\newcommand{\N}{\mathbf{N}}
\newcommand{\RR}{\mathcal R}
\newcommand{\FF}{\mathcal F}
\newcommand{\BB}{\mathcal B}
\newcommand{\ve}{\vee}
\newcommand{\aut}{\mathcal A}
\newcommand{\dqr}{\delta_{\mathfrak q}^{\mathbb R}}

\newcommand{\abs}[1]{\lvert#1\rvert}

\newcommand{\ver}{\phantom{\begin{array}{c}
                         t\\
                         t\\
                         t\\
                         t\\
                         t
                         \end{array}
                        }
                   }

\newcommand{\la}{\langle}
\newcommand{\ra}{\rangle}
\newcommand{\bp}{\bigskip}
\newcommand{\be}{\begin {equation}}
\newcommand{\ee}{\end {equation}}
\newcommand{\bee}{\begin {equation*}}
\newcommand{\eee}{\end {equation*}}
\newcommand{\nt}{\noindent}

\newtheorem{spherical}{Theorem}[section]

\newtheorem{fjcorrespondence}{Definition}[section]
\newtheorem{cdefinition}[fjcorrespondence]{Lemma}
\newtheorem{fjcorrespondenceofd}[fjcorrespondence]{Definition}
\newtheorem{fj}[fjcorrespondence]{Lemma}
\newtheorem{fjsubmersion}[fjcorrespondence]{Lemma}
\newtheorem{doperator}[fjcorrespondence]{Lemma}
\newtheorem{fj2}[fjcorrespondence]{Proposition}
\newtheorem{ksubmersion}[fjcorrespondence]{Lemma}

\newtheorem{jacobian}{Proposition}[section]
\newtheorem{jacobian2}[jacobian]{Lemma}
\newtheorem{jacobian3}[jacobian]{Lemma}
\newtheorem{jacobian4}[jacobian]{Lemma}
\newtheorem{jacobian5}[jacobian]{Lemma}
\newtheorem{spherical2}[jacobian]{Lemma}
\newtheorem{spherical3}[jacobian]{Lemma}
\newtheorem{spherical4}[jacobian]{Lemma}

\newtheorem{dequation}{Theorem}[section]
\newtheorem{dequation2}[dequation]{Lemma}
\newtheorem{dequation3}[dequation]{Lemma}

\newtheorem{bottom}{Lemma}[section]
\newtheorem{bottom2}[bottom]{Lemma}
\newtheorem{bottom3}[bottom]{Proposition}
\newtheorem{bottom4}[bottom]{Definition}
\newtheorem{bottom5}[bottom]{Proposition}

\newtheorem{eandu}{Proposition}[section]
\newtheorem{coef1}[eandu]{Proposition}
\newtheorem{qassociated}[eandu]{Lemma}
\newtheorem{coef4}[eandu]{Theorem}
\newtheorem{coef01}[eandu]{Lemma}
\newtheorem{coef2}[eandu]{Lemma}

\newtheorem{lipositive}{Definition}[section]
\newtheorem{liflensted}[lipositive]{Lemma}
\newtheorem{contragredient}[lipositive]{Lemma}
\newtheorem{tensorproduct}[lipositive]{Lemma}
\newtheorem{invariant}[lipositive]{Proposition}
\newtheorem{branching}[lipositive]{Proposition}
\newtheorem{unitaryinvariant}[lipositive]{Proposition}
\newtheorem{unitarybranching}[lipositive]{Proposition}
\newtheorem{finitedimensional}[lipositive]{Proposition}
\newtheorem{eisensteinli}[lipositive]{Proposition}
\newtheorem{cohomologicallipositive}[lipositive]{Theorem}
\newtheorem{sl2r}[lipositive]{Proposition}
\newtheorem{trivial}[lipositive]{Proposition}
\newtheorem{l1}[lipositive]{Definition}

\title[Eisenstein Integrals]{Eisenstein Integrals for Theta Stable Parabolic Subalgebras}

\author{Sun Binyong}

\thanks{This paper is a part of the author's doctoral dissertation in HKUST, supervised by Jian-Shu Li.}%

\address{Department of Mathematics, Hong Kong University of Science
and Technology}

\email{sunsun@ust.hk}
\date{Nov. 28, 2004}
%
\maketitle
%


\begin{abstract}
Let $G$ be a connected semisimple Lie group with a finite center,
with a maximal compact subgroup $K$. There are two general methods
to construct representations of $G$, namely, ordinary parabolic
induction and cohomological parabolic induction. We define
Eisenstein integrals relative to cohomological inductions which
generalize Flensted-Jensen's fundamental functions for discrete
series. They are analogous to Harish-Chandra's Eisenstein integrals
related to ordinary inductions. We introduce the notion of Li
positivity of a $K$ type in a representation of $G$ which is
extremely useful in the study of branching laws. As an application
of our integral, we show that the minimal $K$ types of many
interesting representations are Li positive. These include all
irreducible unitary representations with nonzero cohomology.
\end{abstract}

\section{Introduction and notations}\label{S:In}

\subsection{Complexifications of groups}

Whenever $H$ is a Lie group, we denote by $H_\C$ the universal
complexification of $H$; $u_H: H\rightarrow H_\C$ the canonical
homomorphism; $\bar \quad : H_\C\rightarrow H_\C$ the
anti-holomorphic automorphism on $H_\C$ which is identity on
$u_H(H)$. Recall that the universal complexification map $u_H$ is
defined by the following property: If $H'$ is a complex Lie group
and $\phi : H\rightarrow H'$ is a Lie group homomorphism, then there
is a unique holomorphic homomorphism $\phi': H_\C \rightarrow H'$
such that $\phi'\circ u_H=\phi$ \cite{Ho}. Notice that although
universal complexifications are only defined for connected Lie
groups in \cite{Ho}, the definitions and the results can be easily
generalized to non-connected Lie groups \cite{Ne}. The followings
are well known:
\begin{enumerate}
              \item
                   The universal complexification of a connected
                   Lie group is connected.
              \item
                   If $H$ is a real reductive group in Harish-Chandra's class, the $u_H$ has a
                   finite kernel and the differential of $u_H$ is a
                   comlexification of the Lie algebra of $H$.
              \item
                   If $H$ is compact, then $u_H$ is injective.

\end{enumerate}

\subsection{Conventions}
If $V$ is a finite dimensional complex vector space and $X$ is a
real analytic manifold, we denote by $\con^\infty(X;V)$ the space of
all smooth functions on $X$ with values in $V$, and
$\con^\omega(X;V)$ for real analytic functions. When $H$ is a Lie
group, $\c$ is a complex subspace of a complex Lie algebra, and $H$
acts on $\c$ by certain adjoint action, we denote by $\delta_\c$ the
non-unitary character on $H$ defined by the determinant of the
action on $\c$. For any finite dimension continuous representation
of a Lie group, say, $\tau$ of $H$, we use the same symbol $\tau$ to
indicate its differential as well as its holomorphic extension to
$H_\C$. Whenever $W$ is a set, we write $1_W$ for the identity map
on $W$, or just $1$ when no confusion is possible. If $H$ is a
compact Lie group, $\alpha$ an $H$ type, i.e., an equivalent class
of irreducible finite dimensional continuous representations of $H$,
and $M$ a $H$ module, we write $M(\alpha)$ for the $\alpha$-isotypic
component of $M$. Finally, we always use the normalized Haar measure
on a compact group for integration.

\subsection{Notations}\label{S:Notations}
The following notations will be fixed throughout this paper unless
otherwise mentioned. Let $G$ be a connected Lie group with Lie
algebra $\g_0$. Let $\g$ be the complexfication of $\g_0$. We assume
$\g_0$ is reductive and the connected Lie subgroup of $G$ with Lie
algebra $[\g_0,\g_0]$ has a finite center. Consequently, this
subgroup is closed in $G$. Let $\theta$ be a Cartan involution on
$G$ and let $K$ be the corresponding maximal compact subgroup of
$G$. $\theta$ induces an automorphism on $\g$ which we still denote
by $\theta$. Let
\[
    \g_0=\k_0+\p_0 \qquad \textrm{and} \qquad \g=\k+\p
\]
be the usual decompositions corresponding to $\theta$. Denote by
\[
  \bar{\theta}=\theta\circ \bar \quad=\bar \quad \circ \theta
\]
a conjugate linear automorphism of $\g$. Notice that the above
"$\bar \quad$" means the complex conjugation of $\g$ with respect to
the real form $\g_0$. $\theta$ induces a holomorphic automorphism on
$G_\C$. We still call it $\theta$. Still denote by
\[
  \bar{\theta}=\theta\circ \bar \quad=\bar \quad \circ \theta
\]
the anti-holomorphic automorphism of $G_\C$. Let $U$ be the subgroup
of $G_\C$ fixed by $\bar{\theta}$. Notice that $U$ is connected; and
$U$ is a maximal compact subgroup of $G_\C$ when $G$ has a compact
center.

Let $\q$ be a parabolic subalgebra of $\g$. We assume $\q$ is real
or $\theta$ stable. These are the two cases which are extremely
interested in representation theory. Here $\q$ is real means $\bar
\q=\q$; $\q$ is $\theta$ stable means $\theta(\q)=\q$ and $\q\cap
\bar \q$ is a Levi factor of $\q$. In both cases we define
\[
     G'=\Norma_G(\q)\cap \Norma_G(\bar{\theta}(\q)).
\]
Here $\Norma_G(\q)$ means the normalizer of $\q$ in $G$, etc. Then
$G'$ is a $\theta$ stable real reductive group in Harish-Chandra's
class. Let $\g_0'$ be the Lie algebra of $G'$ and $\g'$ be the
complexification of $\g_0'$. Then $\g'$ is a Levi factor of $\q$.
Let $\n$ be the nilpotent radical of $\q\cap [\g,\g]$. Then
\[
  \q=\g'\oplus \n.
\]
Denote by $N$ the connected complex subgroup of $G_\C$ with Lie
algebra $\n$. Denote by $K'=K\cap G'$. We note that $G'$ is also
connected in the case of $\theta$ stable $\q$. We define two maps
\[
  \theta,\, \bar{\theta} : G'_\C\rightarrow G'_\C
\]
and a subgroup $U'\subset G'_\C$ analogously.

\subsection{Eisenstein integrals}\label{S:Eint}

When $\q$ is real, we can construct representations of $G$ from
representations of $G'$ via ordinary parabolic inductions. When $\q$
is $\theta$ stable, we can construct representations of $G$ from
representations of $G'$ via cohomological parabolic inductions by
using Zuckerman's functors. Harish-Chandra was able to define his
Eisenstein integral which can be used to express the matrix
coefficient of an ordinary parabolic induced representation.
Harish-Chandra has studied Eisenstein integral extensively
\cite{Ha1}, \cite{Ha2}. We will not talk it in detail in this paper.
In stead, we will concentrate on the other case, namely, when $\q$
is $\theta$ stable. Our integral is a generalization of
Flensted-Jensen's fundamental functions for discrete series.
Harish-Chandra's is a generalization of elementary spherical
functions for spherical principle series. Now we are going to have a
detailed description of our integral, comparing to Harish-Chandra's.

Fix a finite dimensional continuous representation $\tau$ of
$K\times K$ on $\V$. Let
\begin{equation}\label{E:Vprime}
   \V'=\{\,v\in \V\mid \tau(X,Y)v=0 \textrm{ for all } X\in \n\cap
   \k, Y\in \bar \n\cap\k\,\}.
\end{equation}
Define a representation $\tau'$ of $K'\times K'$ on $\V'$ by
\begin{equation}\label{E:tauprime}
   \tau'(k,l)v=\delta_{\bar{\n}\cap\p}(k)\delta_{\n\cap\p}(l)\tau(k,l)v.
\end{equation}
We call $\tau'$ the representation $\q$-associated to $\tau$. Notice
that when $\q$ is real, we have $\V'=\V$ and $\tau'=\tau$.

We define three spaces of $\tau$-spherical functions by
\[
   \begin{array}{lcl}
       \medskip
       \con^\infty(G;\tau)&= & \{\,f \in \con^\infty(G;\V)\mid \\
       \medskip
       &&f(k_1 x k_2^{-1})=\tau(k_1, k_2)f(x) \textrm{ for all } k_1,k_2\in K, x\in G\,\};\\
       \medskip
       \con^\infty(G_\C/U;\tau)&= & \{\,f \in \con^\infty(G_\C/U;\V)\mid\\
       \medskip
       && f(kx)=\tau(k, \bar k)f(x) \textrm{ for all } k\in K_\C, x\in G_\C/U\,\};\\
       \medskip
       \con^\infty(G/K;\tau)&= & \{\,f \in \con^\infty(G/K;\V)\mid\\
       && f(kx)=\tau(k, k)f(x)\textrm{ for all } k\in K, x\in G/K\,\}.
   \end{array}
\]
We will prove that these three spaces can be canonically identified
with each other. We similarly define $\con^\infty(G';\tau')$,
$\con^\infty(G'_\C/U';\tau')$ and $\con^\infty(G'/K';\tau')$.

When $\q$ is real, Langlands decomposition enables us to define a
real analytic map
\[
  \begin{array}{cccc}
       H_\q:& G & \rightarrow & G'\times_{K'}K,\\
       & ngk &\mapsto &\textrm{ the class of }(g,k)
  \end{array}
\]
for all $n\in \exp(\n\cap \g_0)$, $g\in G'$, $k\in K$, where
$G'\times_{K'}K$ means the quotient of $G'\times K$ by $K'$ under
the action indicated by the notation. So $G'\times_{K'}K$ is
essentially the closed subset $G'K$ of $G$. When $\q$ is a $\theta$
stable parabolic, denote by $j_\C:G'_\C\rightarrow G_\C$ the
complexification of the embedding $j:G'\rightarrow G$. Again by
Langlands decomposition we define a real analytic map
\[
   \begin{array}{cccc}
          H_\q : & G_\C/U & \rightarrow & G'_\C/U',\\
            & n j_\C(g) U  & \mapsto & gU'
   \end{array}
\]
for all $n\in N$ and $g\in G'_\C$.

We also need normalization factors. When $\q$ is real, define
\[
   \dqr(g,k)=\abs{\delta_\n(g)}^{1/2}, \qquad (g,k)\in G'\times_{K'}K.
\]
When $\q$ is a $\theta$ stable parabolic, define
\[
   \dqr(xU')=\delta_\n(x\bar{\theta}(x^{-1})), \qquad xU'\in G'_\C/U'.
\]
In both cases, $\dqr$ has positive values everywhere.

When $\q$ is real, $\phi\in\con^\infty(G';\tau')$, extend $\phi$ to
$G'\times_{K'}K$ by
\[
   \phi(g,k)=\tau(1,k^{-1})\phi(g).
\]
Harish-Chandra's Eisenstein integral can be formulated as follows:
\[
   E_\q(\phi:x)=\int_K\tau(k^{-1},1)(\dqr\phi)(H_\q(kx))\,dk\qquad x\in G.
\]
It turns out easily that $E_\q(\phi)\in \con^\infty(G;\tau)$. When
$\q$ is a $\theta$ stable parabolic,
$\phi\in\con^\infty(G'_\C/U';\tau')$, we define our integral by
\[
   E_\q(\phi:x)=\int_K\tau(k^{-1},k^{-1})(\dqr\phi)(H_\q(kx))\,dk\qquad x\in G_\C/U.
\]
We still have $E_\q(\phi)\in \con^\infty(G_\C/U;\tau)$ in our case,
but the proof is not as easy as in Harish-Chandra's case. It is in
fact one of our main results.

\begin{spherical}\label{T:spherical}
       When $\q$ is a $\theta$ stable parabolic,
       $E_\q(\phi)\in \con^\infty(G_\C/U;\tau)$ for all
       $\phi\in\con^\infty(G'_\C/U';\tau')$.
\end{spherical}

We will prove the above theorem later.

Denote by $\tau|_K$ the representation of $K$ on $\V$ given by
\[
   \tau|_K(k)=\tau(k,k), \qquad k\in K.
\]
If $k\in K$, denote by
\[
  \begin{array}{l}
     T_k(x)=kxk^{-1}, \quad \textrm{if}\quad x\in G, \quad
      \textrm{and}\\
     T_k(xU)=kxU, \quad \textrm{if}\quad xU\in G_\C/U.
  \end{array}
\]

Then it's clear that Harish-Chandra's integral and ours can be
formulated in the same appearance:
\[
   E_\q(\phi:x)=\int_K\tau|_K(k^{-1})((\dqr\phi)\circ H_\q\circ T_k)(x)\,dk.
\]
In both cases, the normalization $\dqr$ is used to make the
integral compatible with the usual normalization of parabolic
inductions. Matrix coefficients of parabolic induced
representations can be expressed as Harish-Chandra's Eisenstein
integrals. We will prove that the matrix coefficients of the
bottom layers of cohomologically induced representations can be
expressed as our Eisenstein integral. Combining these two, we can
find an integral representation of the matrix coefficient of every
minimal $K$ type in an arbitrary irreducible $(\g,K)$ module.

The author would like to thank Prof. Jian-Shu Li, for his
initiating of this project and his guidance on study.

\section{Flensted-Jensen's duality}
The purpose of this section is to understand the correspondence
between $\con^\infty(G;\tau)$ and $\con^\infty(G_\C/U;\tau)$
mentioned in last section. The ideas come from \cite{FJ}.
\subsection{Complexification squares}

Denote by $\S_A$ a commutative diagram
\[
   \begin{CD}
          A_1          @>>>                     A_2 \\
          @VVV                               @VVV\\
          A_3       @>>>                    A_4.
   \end{CD}
\]
We call $\S_A$ a complexification square if the followings are
satisfied:
\begin{enumerate}
              \item
                   $A_1$, $A_2$, $A_3$ are nonempty connected real
                   analytic manifolds, $A_4$ is a connected complex manifold;
              \item
                   All the arrows in $\S_A$ are real analytic
                   maps;
              \item
                   The differential of the map $A_2\rightarrow A_4$
                   is a complexification of real vector space everywhere;
             \item The same for the map $A_3\rightarrow A_4$.
\end{enumerate}
The last two conditions essential say that locally $A_2$ and $A_3$
are totally real submanifolds of $A_4$.

We say $\S_A$ is a complexification square of groups if
\begin{enumerate}
              \item
                   $A_1$, $A_2$, $A_3$ are connected Lie groups;
                   $A_4$ is a connected complex Lie group.
              \item
                   All the arrows in $\S_A$ are Lie group homomorphisms;
              \item
                   The differential of the homomorphism $A_2\rightarrow A_4$
                   is a complexification of the Lie algebra of
                   $A_2$;
              \item
                   The same for the map $A_3\rightarrow A_4$.
\end{enumerate}
Therefore a complexification square of groups is a complexification
square.

Let $\S_A$ be a complexification square (complexification square of
groups) as before. Let $\S_B$
\[
   \begin{CD}
          B_1          @>>>                     B_2 \\
          @VVV                               @VVV\\
          B_3       @>>>                     B_4.
   \end{CD}
\]
be another complexification square (complexification square of
groups ). We define the product $\S_A\times \S_B$ to be the
commutative diagram
\[
   \begin{CD}
          A_1\times B_1          @>>>                     A_2\times B_2 \\
          @VVV                               @VVV\\
          A_3\times B_3       @>>>                    A_4\times B_4,
   \end{CD}
\]
where the arrows are the corresponding products. Then $\S_A\times
\S_B$ is still a complexification square (complexification square of
groups). We define an homomorphism from $\S_A$ to $\S_B$ to be four
maps
\[
  f_i: A_i\rightarrow B_i\qquad i=1,2,3,4,
\]
such that
\begin{enumerate}
        \item
            $f_1$, $f_2$, $f_3$ are real analytic maps (Lie group homomorphisms),
       \item
            $f_4$ is a complex analytic map (Complex analytic group homomorphism),
       \item
            the diagram
    \[
     \begin{CD}
          A_1      @>>>                        A_2\\
          \left \downarrow\ver\right.&
                \begin{array}{rcl}
                       \searrow f_1 &\phantom{\rightarrow}             &\swarrow f_2\\
                       \quad B_1    &\rightarrow                         &B_2\quad  \\
                       \quad\downarrow &                                 &\downarrow\\
                       \quad B_3 &\rightarrow                         &B_4\\
                       \nearrow f_3  &\phantom{\rightarrow}             &\nwarrow f_4
                \end{array}&
         \left.\ver \right \downarrow\\
         A_3 @>>>                         A_4.
    \end{CD}
   \]
  commutes.
\end{enumerate}
Denote by
  \[
    f=(f_1,f_2,f_3,f_4): \S_A\rightarrow\S_B
  \]
for the homomorphism.

Assume $\S_A$ is a complexification square of groups, $\S_B$ is a
complexification square, and
 \[
    f=(f_1,f_2,f_3,f_4): \S_A\times \S_B \rightarrow\S_B
  \]
is a homomorphism of complexification squares, we say $f$ is an
analytic action of $\S_A$ on $\S_B$ if all $f_i$, $i=1,2,3,4$ are
group actions. We define the concept of a complexification square of
Lie algebras in the obvious way. If $\S_A$ is a complexification
square of groups, take the Lie algebra of each component of $\S_A$,
we get a complexification square of Lie algebras. We call it the Lie
algebra of $\S_A$.

\subsection{Correspondence of functions}

Assume $\S_A$
\[
   \begin{CD}
          A_1          @>u_2>>                     A_2 \\
          @Vu_3VV                               @VVv_2V\\
          A_3       @>>v_3>                    A_4.
   \end{CD}
\]
is a complexification square, and $V$ is a finite dimensional
complex vector space. Fix a base point $a_1\in A_1$. Denote by
\[
  a_2=u_2(a_1),\quad a_3=u_3(a_1),\quad a_4=v_2(a_2).
\]

\begin{fjcorrespondence}\label{T:fjcorrespondence}
Let $\phi\in \con^\omega(A_2;V)$, and $\psi\in\con^\omega(A_3;V)$.
$\phi$ and $\psi$ are said to correspond to each other under $\S_A$
(or $\phi\equiv_{\S_A}\psi$ in notation) if the following condition
is satisfied:

There are open neighborhoods $U_i$ of $a_i$ in $A_i$, $i=2,3,4$ and
a holomorphic map $f$ from $U_4$ to $V$ such that
\begin{enumerate}
       \item
            $v_2(U_2)\subset U_4$, $v_3(U_3)\subset U_4$, and
       \item
            $f\circ (v_2|_{U_2})=\phi|_{U_2}$, $f\circ (v_3|_{U_3})=\psi|_{U_3}$.
\end{enumerate}
\end{fjcorrespondence}
Notice that any function in $\con^\omega(A_2;V)$ corresponds to at
most one function in $\con^\omega(A_3;V)$, and vice versa.

\begin{cdefinition}
The above definition is independent of the choice of the base
point $a_1\in A_1$.
\end{cdefinition}

\proof Let $\phi\in \con^\omega(A_2;V)$, and
$\psi\in\con^\omega(A_3;V)$. In this proof, we use
\[
  \phi\equiv_{a_1}\psi
\]
to indicate that $\phi$ and $\psi$ correspond to each other under
the above definition. Denote by
\[
  A_0=\{\,b\in A_1\mid \phi\equiv_{b}\psi \,\}.
\]
We only need to show that $A_0$ is either $A_1$ or the empty set.
It's clear that $A_0$ is open in $A_1$. We should prove that $A_0$
is closed. Let $b_1$ be a point in the closure of $A_0$ in $A_1$.
Denote by
\[
  b_2=u_2(b_1),\quad b_3=u_3(b_1),\quad b_4=v_2(b_2).
\]
Choose open connected neighborhoods $W_i$ of $b_i$ in $A_i$,
$i=1,2,3,4$, so that
\[
    \begin{array}{l}
     u_2(W_1)\subset W_2,\qquad u_3(W_1)\subset W_3, \quad \textrm{and}\\
     v_2(W_2)\subset W_4,\qquad v_3(W_3)\subset W_4.
    \end{array}
\]
We choose them small enough, then there are holomorphic functions
$g$ and $h$ on $W_4$ satisfies:
\[
  \begin{array}{l}
   g\circ (v_2|_{W_2})=\phi|_{W_2}, \qquad \textrm{and}\\
   h\circ (v_3|_{W_3})=\psi|_{W_3}.
  \end{array}
\]
Since $b_1$ is in the closure of $A_0$, there is a point $a_1\in
W_1\cap A_0$. Denote by
\[
  a_2=u_2(a_1),\quad a_3=u_3(a_1),\quad a_4=v_2(a_2).
\]
By definition, there are open connected neighborhoods $U_i$ of $a_i$
in $W_i$, $i=1,2,3,4$, so that there is a holomorphic function $f$
on $U_4$ satisfies:
\begin{enumerate}
       \item
            $v_2(U_2)\subset U_4$, $v_3(U_3)\subset U_4$, and
       \item
            $f\circ (v_2|_{U_2})=\phi|_{U_2}$, $f\circ (v_3|_{U_3})=\psi|_{U_3}$.
\end{enumerate}
Then we have
\[
  f\circ (v_2|_{U_2})=g|_{U_4}\circ (v_2|_{U_2}).
\]
Therefore
\[
  f=g|_{U_4}.
\]
Similarly
\[
  f=h|_{U_4}.
\]
Now we conclude that $g=h$ and therefore $b_1\in A_0$. This proves
that $A_0$ is closed.

\qed

The following lemma will be used later, we omit its easy proof.

\begin{fjsubmersion}\label{T:fjsubmersion}
Assume both $\S_A$ and $\S_B$ are complexification squares. Let
\[
   f=(f_1,f_2,f_3,f_4): \S_A\rightarrow\S_B
\]
be a homomorphism from $\S_A$ to $\S_B$. Assume that both $f_2$ and
$f_3$ are submersions. Let $\phi\in \con^\omega(B_2;V)$ and
$\psi\in\con^\omega(B_3;V)$. Then
\[
  \phi\equiv_{\S_B}\psi\quad \textrm{if and only if}\quad \phi\circ f_2\equiv_{\S_A}\psi\circ
  f_3.
\]
\end{fjsubmersion}

\subsection{Differential operators}
We can generalize the above considerations to differential
operators. Let $A$ be a real analytic manifold, $B$ a complex
manifold, and
\[
  u : A\rightarrow B
\]
a real analytic map whose differential is a complexification of real
vector space everywhere. Let $D$ be a holomorphic differential
operator on $B$. Then there is a unique analytic differential
operator $D'$ on $A$ satisfies the following:

For every open subset $U_A$ of $A$, $U_B$ of $B$, and every
holomorphic function $f$ on $U_B$, if $u(U_A)\subset U_B$, then
\[
  D'(f\circ u|_{U_A})=D(f)\circ u|_{U_A}.
\]
Denote by $u^*(D)$ the differential operator $D'$.

Assume $\S_A$
\[
   \begin{CD}
          A_1          @>u_2>>                     A_2 \\
          @Vu_3VV                               @VVv_2V\\
          A_3       @>>v_3>                    A_4.
   \end{CD}
\]
is a complexification square as before. Fix a base point $a_1\in
A_1$. Recall that
\[
  a_2=u_2(a_1),\quad a_3=u_3(a_1),\quad a_4=v_2(a_2).
\]

\begin{fjcorrespondenceofd}\label{T:fjcorrespondenceofd}
Let $D_2$ be an analytic differential operator on $A_2$, and $D_3$
an analytic differential operator on $A_3$. $D_2$ and $D_3$ are said
to correspond to each other under $\S_A$ (or $D_2\equiv_{\S_A}D_3$
in notation) if the following condition is satisfied:

There are open neighborhoods $U_i$ of $a_i$ in $A_i$, $i=2,3,4$ and
a holomorphic differential operator $D$ on $U_4$ such that
\begin{enumerate}
       \item
            $v_2(U_2)\subset U_4$, $v_3(U_3)\subset U_4$, and
       \item
            $(v_2|_{U_2})^*(D)=D_2|_{U_2}$, $(v_3|_{U_3})^*(D)=D_3|_{U_3}$.
\end{enumerate}
\end{fjcorrespondenceofd}

Again this definition is independent of the choice of the base point
$a_1\in A_1$. An analytic differential operator on $A_2$ correspond
to at most one analytic differential operator on $A_3$, and vice
versa.

Let $V$ be a finite dimensional complex vector space. The following
lemma is obvious.

\begin{fj}\label{T:fj}
       Let $\phi\in \con^\omega(A_2;V)$, $\psi\in\con^\omega(A_3;V)$.
       Let $D_2$ be an analytic differential operator on $A_2$, and $D_3$
an analytic differential operator on $A_3$. If
\[
  \phi\equiv_{\S_A}\psi, \quad \textrm{and} \quad
  D_2\equiv_{\S_A}D_3,
\]
then
       \[
         D_2(\phi)\equiv_{\S_A}D_3(\psi).
       \]
\end{fj}

\subsection{Some complexification squares}

We have a commutative diagram which we denote by $\S$:
\[
   \begin{CD}
      G/K         @>p_K>>     G  \\
      @Vv_GVV                 @VVu_GV\\
      G_\C/U      @>p_U>>     G_\C,
   \end{CD}
\]
where $v_G$ is the map induced by $u_G$; $p_K$ is defined by
\[
  p_K(xK)=x\theta(x^{-1}), \quad x\in G;
\]
$p_U$ is defined by
\[
  p_U(xU)=x\bar{\theta}(x^{-1}), \quad x\in G_\C.
\]
$\S$ is clearly a complexification square. Denote by $\S_G$ the
following complexification square of groups:
\[
   \begin{CD}
      G         @>1\times\theta>>     G\times G  \\
      @Vu_GVV                         @VVu_G\times u_GV\\
      G_\C      @>1\times\bar{\theta}>>     G_\C\times G_\C.
   \end{CD}
\]
We have a natural analytic action $T$ of $\S_G$ on $\S$. The action
of $G\times G$ on $G$ is given by
\[
     T_{g,l}x=gxl^{-1} \qquad g,l,x\in G.
\]
The action of $G_\C\times G_\C$ on $G_\C$ is defined similarly. The
other two actions are defined by the natural left translations.

Denote by $\S_\g$ the Lie algebra of $\S_G$, it is the following
complexification square of Lie algebras:
\[
   \begin{CD}
      \g_0         @>1\times\theta>>     \g_0\times \g_0  \\
      @VVV                               @VVV\\
      \g      @>1\times\bar{\theta}>>          \g\times \g.
   \end{CD}
\]
Let $V$ be a finite dimensional complex vector space as before. The
action $T$ of $\S_G$ on $\S$ induces a smooth representation $T$ of
$G\times G$ on $\con^\infty(G;V)$ by
\[
  \begin{array}{rl}
   (T_{g,h}\phi)(x)=\phi(T_{g^{-1},h^{-1}}(x))&=\phi(g^{-1}xh),\\
    &g,h,x\in G,\,\phi\in\con^\infty(G;V)
   \end{array}
\]
and also a smooth representation $T$ of $G_\C$ on
$\con^\infty(G_\C/U;V)$ by
\[
  \begin{array}{rl}
   (T_g\psi)(x)=\psi(T_{g^{-1}}(x))&=\psi(g^{-1}x),\\
    &g\in G_\C,\,x\in
   G_\C/U,\,\psi\in\con^\infty(G_\C/U;V).
  \end{array}
\]
Take differential of these two representations, and use the
complexifications in the diagram $\S_\g$, we have an action $T$ of
$\U(\g)\otimes \U(\g)$ on $\con^\infty(G;V)$ and
$\con^\infty(G_\C/U;V)$. Notice that here we use the canonical
identification of $\U(\g)\otimes \U(\g)$ with $\U(\g\times \g)$.
This $T$ is in fact the action of analytic differential operators.

\begin{doperator}\label{T:doperator}
Let $X\otimes Y\in \U(\g)\otimes \U(\g)$. Let $D_G$ be the
differential operator on $G$ determined by $T_{X\otimes Y}$. Let
$D_{G_\C/U}$ be the differential operator on $G_\C/U$ determined by
$T_{X\otimes Y}$. Then
\[
  D_G\equiv_{\S} D_{G_\C/U}.
\]
\end{doperator}

\proof Let $D_\C$ be the holomorphic differential operator on $G_\C$
determined by $X\otimes Y$ and the holomorphic action $T$ of $G_\C
\times G_\C$ on $G_\C$. Then
\[
  (u_G)^*(D_\C)=D_G, \quad \textrm{and} \quad (p_U)^*(D_\C)=D_{G_\C/U}.
\]
\qed

\subsection{Flensted-Jensen's duality for groups}
Let $\tau$ be a finite dimensional continuous representation of
$K\times K$ on $\V$ as before. $\con^\omega(G;\tau)$,
$\con^\omega(G_\C/U;\tau)$ and $\con^\omega(G/K;\tau)$, etc, are
defined in an obvious way.

The following is the main result of this section.

\begin{fj2}\label{T:fj2}
       \begin{enumerate}
              \item
                   The pull back of $p_K$ and $v_G$ induces canonical isomorphisms of vector spaces:
                   \[
                      \con^\infty(G;\tau)=\con^\infty(G/K;\tau)=\con^\infty(G_\C/U;\tau).
                   \]
              \item
                   \[
                      \con^\omega(G;\tau)=\con^\omega(G/K;\tau)=\con^\omega(G_\C/U;\tau)
                   \]
                   under the above identifications.
              \item
                  If $\phi\in \con^\omega(G;\tau)$ and $\psi\in\con^\omega(G_\C/U;\tau)$, then
                  \[
                     \phi\equiv_{\mathbf S}\psi \quad \textrm{if and only if}
                     \quad \phi\circ p_K=\psi\circ v_G.
                  \]
       \end{enumerate}
\end{fj2}

Notice that
\[
  \con^\omega(G;\tau)=\con^\omega(G_\C/U;\tau)
\]
is also an identification of $\U(\g)^K\otimes \U(\g)^K$ modules.

\subsection{More complexification squares}

Denote by $\S_K$  the following complexification square of groups
\[
   \begin{CD}
          K          @>1\times 1>>              K\times K  \\
          @Vu_KVV                               @VVu_K\times u_KV\\
          K_\C       @>1\times \bar{\quad}>>    K_\C\times K_\C.
   \end{CD}
\]
Denote by $\S_{G/K}$ the complexification square
\[
   \begin{CD}
          G/K          @>1>>                  G/K \\
          @V1VV                               @VVV\\
          G/K          @>>>                   G_\C/K_\C,
   \end{CD}
\]
where the right vertical arrow and the bottom horizontal arrow are
the maps induced by $u_G$. We define four maps
\[
   \begin{array}{llll}
          \pi_1 : & K\times G/K & \rightarrow & G/K,  \\
           \medskip
                        & (k,xK)      & \mapsto     & kxK;  \\

          \pi_2 : & (K\times K)\times G/K & \rightarrow & G,  \\
           \medskip             & (k_1,k_2,xK)      & \mapsto     & k_1 x\theta(x^{-1})k_2^{-1};  \\

          \pi_3 : & K_\mathbb{C}\times G/K & \rightarrow & G_\mathbb{C}/U,  \\
          \medskip              & (k,xK)                 & \mapsto     & ku_G(x)U;  \\

          \pi_4 : & (K_\C\times K_\C)\times G_\C/K_\C & \rightarrow & G_\C,  \\
                        & (k_1,k_2,xK_\C)      & \mapsto     & k_1 x\theta(x^{-1})k_2^{-1}.  \\
   \end{array}
\]
Here $\pi_1$, $\pi_2$, $\pi_3$ are real analytic; $\pi_4$ is
holomorphic.

It's routine to check that the following diagram commutes:
\[
   \begin{CD}
          K\times G/K   @>>>                        (K\times K)\times G/K\\
          \left \downarrow\ver\right.&
                \begin{array}{rcl}
                       \searrow \pi_1 &\phantom{\rightarrow}             &\swarrow\pi_2\\
                       \quad G/K    &\rightarrow                         &G\quad  \\
                       \quad\downarrow &                                 &\downarrow\\
                       \quad G_\C/U &\rightarrow                         &G_\C\\
                       \nearrow\pi_3  &\phantom{\rightarrow}             &\nwarrow\pi_4
                \end{array}&
         \left.\ver \right \downarrow\\
         K_\C\times G/K @>>>                         (K_\C\times K_\C)\times G_\C/K_\C.
  \end{CD}
\]
Here the outer side square is the product $\S_K\times \S_{G/K}$, and
the inner square is $\S$. Therefore by putting all the four $\pi$'s
together, we have a homomorphism
\[
  \pi=(\pi_1,\pi_2,\pi_3,\pi_4) : \S_K\times \S_{G/K}\rightarrow \S.
\]

Denote by $\S_0$ the following complexification square of groups
\[
   \begin{CD}
          K          @>>>              K  \\
          @VVV                               @VVu_KV\\
          K          @>u_K>>          K_\C.
   \end{CD}
\]
If $a\in K$, we define the actions of $a$ on $K\times G/K$,
$(K\times K)\times G/K$, and $K_\mathbb{C}\times G/K$ by
\[
  \begin{array}{lll}
                    (k,xK)      & \mapsto     & (ka^{-1},axK),  \\
                    (k_1,k_2,xK)      & \mapsto     &       (k_1 a^{-1},k_2 a^{-1},axK),\textrm{ and }\\
                    (k,xK)      & \mapsto     & (ka^{-1},axK),

  \end{array}
\]
respectively. If $a\in K_\C$, we define the actions of $a$ on $(K_\C
\times K_\C)\times G_\C/K_\C$ by
\[
  \begin{array}{lll}
                    (k_1,k_2,xK_\C)      & \mapsto     &       (k_1 a^{-1},k_2
                    a^{-1},axK).
  \end{array}
\]
It's again routine to check that we have an analytic action of
$\S_0$ on $\S_K\times \S_{G/K}$ by the above actions.

The following lemma says that $\S$ is the quotient of $\S_K\times
\S_{G/K}$ by $\S_0$, except the complex component. It is implicit in
\cite{FJ}. We omit its proof.

\begin{ksubmersion}\label{T:ksubmersion}
The maps $\pi_1$, $\pi_2$, and $\pi_3$ are surjective real analytic
submersions. Each fibre of theses maps is a $K$ obit of the actions
defined above.
\end{ksubmersion}

\subsection{Proof of Proposition \ref{T:fj2}}

Denote by
\begin{eqnarray*}
        &&\quad\con^\infty((K\times K)\times G/K;\tau)\\
        &&=\{\,\phi\in \con^\infty((K\times K)\times
        G/K;\V)\mid \phi(k_1,k_2,x)=\tau(k_1,k_2)\phi(1,1,x),\\
        && \qquad\phi(k_1a^{-1},k_2a^{-1},ax)=\phi(k_1,k_2,x),
        \quad k_1,k_2,a\in K,\,x\in G/K\,\}.
\end{eqnarray*}
By using the pull back of $\pi_2$ we have an identification
\[
  \con^\infty(G;\tau)=\con^\infty((K\times K)\times G/K;\tau).
\]
We identify $G/K$ with the subset $\{1\}\times\{1\}\times G/K$ of
$(K\times K)\times G/K$. Every function in $\con^\infty((K\times
K)\times G/K;\tau)$ is determined by its restriction to $G/K$. It's
easy to see that by using the restriction we have an identification
\[
  \con^\infty((K\times K)\times G/K;\tau)=\con^\infty(G/K;\tau).
\]

Similarly denote by
\begin{eqnarray*}
        &&\quad\con^\infty(K_\C\times G/K;\tau)\\
        &&=\{\,\phi\in \con^\infty(K_\C\times G/K;\V)\mid \phi(k,x)
          =\tau(k,\bar k)\phi(1,x),\\
        &&\qquad\phi(ka^{-1},ax)=\phi(k,x), \qquad k\in K_\C,\,a\in K,\,x\in G/K\,\}.
\end{eqnarray*}
Then by using the pull back of $\pi_3$ and a restriction map we have
canonical identifications
\[
 \con^\infty(G_\C/U;\tau)=\con^\infty(K_\C\times
 G/K;\tau)=\con^\infty(G/K;\tau).
\]
Keep these identifications in mind, we get the first assertion of
Proposition \ref{T:fj2}. It's clear that under these
identifications, real analytic functions correspond to real
analytic functions. This proves the second assertion.

For the third assertion of Proposition \ref{T:fj2}. Let $\phi\in
\con^\omega(G;\V)$ and $\psi\in\con^\omega(G_\C/U;\V)$. If
$\phi\equiv_{\mathbf S}\psi $, then it's clear that
\[
   \phi\circ p_K=\psi\circ v_G.
\]
Now assume
\[
  \phi\circ p_K=\psi\circ\ v_G=f_0.
\]
Choose a connected open neighborhood $C_0$ of $1K$ in $G/K$ and a
connected open neighborhood $C$ of $1K_\C$ in $G_\C/K_\C$ so that
$f_0|_{C_0}$ extends to a holomorphic function $f$ on $C$. Extend
$f$ to $K_\C\times K_\C\times C$ by
\[
  f(k_1,k_2,x)=\tau(k_1,k_2)f(x).
\]
By this extended $f$, we easily get that $\phi\circ \pi_2$ and
$\psi\circ \pi_3$ correspond to each other under $\S_K\times
\S_{G/K}$. By Lemma \ref{T:fjsubmersion} and Lemma
\ref{T:ksubmersion} , we know $\phi\equiv_{\mathbf S}\psi$. This
finishes the proof of the third assertion of Proposition
\ref{T:fj2}.

\section{A proof of Theorem \ref{T:spherical}}

\subsection{An integral formula}
In this subsection, we establish an integral formula which is
crucial for the proof of Theorem \ref{T:spherical}. The formula is a
stronger version of Corollary 11.40 of \cite{KV}. We only need to
apply the formula to the complex group $K_\C$. We formulate it in
its full generality as it is interesting in itself. The notations in
this subsection are not used in other parts of this paper.

Let $G$ be a real reductive Lie group in Harish-Chandra's class,
$\theta$ a Cartan involution on $G$. Let $\q_0$ be a parabolic
subalgebra of $\g_0=\Lie(G)$. We introduce the following notations:
\[
   \begin{array}{l}
          K \textrm{ is the maximal compact subgroup of $G$ fixed by } \theta.\\
          \g_0=\k_0+\p_0 \textrm{ is the Cartan decomposition corresponds to }\theta\\
          Q \textrm{ is the normalizer of }\q_0 \textrm{ in } G.\\
          Q=LN \textrm{ is the Levi decomposition so that $L$ is $\theta$ stable.}.\\
          \q_0=\l_0+\n_0 \textrm{ is the corresponding Levi decomposition at the Lie algebra level.}\\
          \l_0'=\l_0\cap \p_0. \\
          \delta:L\rightarrow \R^\times\textrm{ is the homomorphism defined by the absolute value
          of }\\
          \qquad\textrm{the determinant of the adjoint representation of $L$ on $\n_0$.}
   \end{array}
\]

By Langlands decomposition, the map
\begin{equation}\label{E:Langlands}
   \begin{array}{rcl}
         N\times\l_0'\times K & \rightarrow & G,\\
         (n,X,k) & \mapsto &  n\exp(X)k
   \end{array}
\end{equation}
is a real analytic diffeomorphism. We define two maps $\kappa :
G\rightarrow K$, and $I: G\rightarrow L$ by $\kappa(nak)=k,
I(nak)=a$ for all $n\in\ N$, $a\in\exp(\l_0')$, $k\in K$. If $a\in
G$, we define a map $\kappa^a : K\rightarrow K$ by
$\kappa^a(k)=\kappa(ka)$.

The purpose of this subsection is to prove

\begin{jacobian}\label{T:jacobian}
If $f$ is a continuous function on $K$, then
\[
   \int_K f(\kappa^a(k))\,dk=\int_K
   f(k)\delta(I(ka^{-1}))\,dk.
\]
\end{jacobian}

We need some lemmas.

\begin{jacobian2}
$\kappa^a$ and $\kappa^{a^{-1}}$ are inverse to each other. Hence
$\kappa^a$ is an analytic diffeomorphism.
\end{jacobian2}

\proof Let $k\in K$ and $ka=n_1s_1k_1$, where $n_1\in N,
s_1\in\exp(\l_0'),$ and $k_1\in K$. Then we have
\[
   k_1a^{-1}=s_1^{-1}n_1^{-1}k=(s_1^{-1}n_1^{-1}s_1)s_1^{-1}k.
\]
Hence
\[
   \kappa^{a^{-1}}\kappa^a(k)=\kappa^{a^{-1}}(k_1)=\kappa(k_1a^{-1})=k.
\]
Change $a$ to $a^{-1}$, we get $\kappa^a\kappa^{a^{-1}}(k)=k$.
\qed

If $a\in G$, we denote by $T_a$ the right translation of $a$ on
$Q\backslash G$. We give $Q\backslash G$ the unique $K$ invariant
measure with total mass $1$.

\begin{jacobian3}\label{T:jacobian3}
If $k\in K$, then the Jacobian of $T_a$ at $Qk$ is
\[
   J_{T_a}(Qk)=\delta(I(ka)).
\]
\end{jacobian3}

This should be known, we give a proof for the sake of completeness.

\proof Let $ka=ns\kappa(ka)$, where $n\in N$ and $s=I(ka)$. We have
a commutative diagram
\[
   \begin{CD}
          Q\backslash G     @>T_a>> Q\backslash G  \\
          @VT_{k^{-1}}VV    @VVT_{\kappa(ka)^{-1}}V\\
          Q\backslash G     @>T_{ns}>> Q\backslash G.
   \end{CD}
\]
Since $T_{k_{-1}}$ and $T_{\kappa(ka)^{-1}}$ preserve the measure,
\[
   J_{T_a}(Qk)=J_{T_{ns}}(Q1).
\]
The right hand side is the absolute value of the determinant of the
tangent map $\dif T_{ns}|_Q$ at the point $Q$. We have another
commutative diagram
\[
   \begin{CD}
          G               @>\Ad_{(ns)^{-1}}>>    G  \\
          @VVV            @VVV \\
          Q\backslash G   @>T_{ns}>>             Q\backslash G.
   \end{CD}
\]
Take the tangent map at the identity, we get
\[
   \begin{CD}
          \g_0           @>\Ad_{(ns)^{-1}}>>             \g_0 \\
          @VVV                                           @VVV \\
          \g_0/\q_0      @>\dif\!T_{ns}|_{Q1}>>   \g_0/\q_0.
   \end{CD}
\]
Hence
\[
   J_{T_a}(Qk)=J_{T_{ns}}(Q)=\abs{\det(\dif T_{ns}|_Q)}
   =\frac{\abs{\det(\Ad_{(ns)^{-1}})}}{\abs{\det(\Ad_{(ns)^{-1}}|\q_0)}}
   =\frac{1}{\delta(s^{-1})}=\delta(s).
\]
\qed

For any $a\in G$, let $J_a$ be the Jacobian of $\kappa^a$.

\begin{jacobian4}
$J_a$ is left $K\cap Q$ invariant.
\end{jacobian4}

\proof
For any $k\in K$, we use $L_k$ to denote the left
translation on $K$ by $k$. The lemma comes from the fact that if
$k\in K\cap Q$, then the following diagram commutes:
\[
   \begin{CD}
          K     @>\kappa^a>>    K  \\
          @VL_kVV               @VVL_kV \\
          K     @>\kappa^a>>    K.
   \end{CD}
\]
We omit the easy proof of this fact.
\qed

\begin{jacobian5}
$J_a(k)=\delta(I(ka))$ for all $k\in K$.
\end{jacobian5}

\proof
Let $\pi$ be the map $K\rightarrow Q\backslash G$,
$k\mapsto Qk$. By the above lemma, there is a function $J'_a$ on
$Q\backslash G$ such that $J'_a\circ\pi=J_a$. For all continuous
function $f$ on $Q\backslash G$, we have
\begin{eqnarray*}
         &&\quad\int_{Q\backslash G}f(T_a(x))J'_a(x)\,dx\\
         &&=\int_K f(T_a(\pi(k)))J'_a(\pi(k))\,dk\\
         &&=\int_K f((\pi(\kappa^a(k))))J_a(k)\,dk  \qquad (T_a\circ\pi=\pi\circ\kappa^a)\\
         &&=\int_K f(\pi(k))\,dk\\
         &&=\int_{Q\backslash G}f(x)\,dx\\
         &&=\int_{Q\backslash G}f(T_a(x))J_{T_a}(x)\,dx.
\end{eqnarray*}
Hence $J'_a=J_{T_a}$. We conclude the proof by Lemma
\ref{T:jacobian3}.
\qed

\it{Proof of Proposition \ref{T:jacobian}.}
\begin{eqnarray*}
          &&\quad\int_K f(\kappa^a(k))\,dk\\
          &&=\int_K f(\kappa^a(\kappa^{a^{-1}}(k)))J_{a^{-1}}(k)\,dk\\
          &&=\int_K f(k)\delta(I(ka^{-1}))\,dk.
\end{eqnarray*}
\qed

\subsection{The proof}

We now return to use the notations in section \ref{S:Eint}. Assume
$\q$ is a $\theta$ stable parabolic. We apply the results obtained
in the last subsection to the group $K_\C$. Let $\k_0'$ be the Lie
algebra of $K'$. Let $N_c$ be the connected subgroup of $K_\C$
with Lie algebra $\n\cap\k$. Now the parabolic subgroup is
$N_cK'_\C$. The map \eqref{E:Langlands} in this case is
\[
   \begin{array}{rcl}
         N_c\times\sqrt{-1}\k_0'\times K & \rightarrow & K_\C,\\
         (n,X,k) & \mapsto &  n\exp(X)k.
   \end{array}
\]
Define $\kappa : K_\C\rightarrow K$; $I: K_\C\rightarrow K'_\C$;
$\kappa^a:K\rightarrow K,\,a\in K_\C$; and $\delta:K'_\C\rightarrow
\R^\times$ as in the last subsection.

\begin{spherical2}\label{T:spherical2}
\[
   H_\q(kax)=(I(\kappa^a(k)a^{-1}))^{-1}H_\q(\kappa^a(k)x)
\]
for all $k\in K, a\in K_\C$ and $x\in G_\C/U$.
\end{spherical2}

\proof Let $ka=n'a'k'$, where $n'\in N_c$, $a'\in
\exp(\sqrt{-1}\k'_0)$, $k'\in K'$. Then
\[
   \kappa^a(k)a^{-1}=k'a^{-1}=a'^{-1}n'^{-1}k=(a'^{-1}n'^{-1}a')a'^{-1}k.
\]
Therefore $I(\kappa^a(k)a^{-1})=a'^{-1}$. As $H_\q$ is $G'_\C$
equivariant and $N$ invariant, we have
 \begin{eqnarray*}
         &&\quad(I(\kappa^a(k)a^{-1}))^{-1}H_\q(\kappa^a(k)x)\\
         &&=a'H_\q(\kappa^a(k)x)=H_\q(a'\kappa^a(k)x)\\
         &&=H_\q(a'k'x)=H_\q(n'^{-1}kax)=H_\q(kax)
 \end{eqnarray*}
\qed

Define the representation $\tau|_{K_\C}$ of $K_\C$ on $\V$ by
$\tau|_{K_\C}(k)=\tau(k,\bar k)$.

\begin{spherical3}\label{T:spherical3}
\[
   (\dqr\phi)(kx)=\delta_{\n\cap\k}(k)\delta_{\bar \n\cap \k}(\bar
          k)\tau|_{K_\C}(k)(\dqr\phi)(x)
\]
for all $k\in K'_\C$, $x\in G'_\C/U'$, and $\phi\in
\con^\infty(G'_\C/U';\tau')$.
\end{spherical3}

\proof
   \begin{eqnarray*}
          &&\quad\dqr\phi(kx)\\
          &&=\dqr(kx)\phi(kx)\\
          &&=\delta_\n(k)\delta_\n({\bar k}^{-1})\dqr(x)\tau'(k,\bar k)\phi(x)\\
          &&=\delta_\n(k)\delta_\n({\bar k}^{-1})\delta_{\bar \n\cap\p}(k)
          \delta_{\n\cap\p}(\bar k)\dqr(x)\tau(k,\bar k)\phi(x)\\
          &&=\delta_{\n\cap\k}(k)\delta_{\bar \n\cap \k}(\bar
          k)\tau|_{K_\C}(k)(\dqr\phi)(x).
  \end{eqnarray*}
\qed

We omit the proof of the following elementary lemma.

\begin{spherical4}\label{T:spherical4}
\[
  \delta(k^{-1})\delta_{\n\cap\k}(k)\delta_{\bar \n\cap \k}(\bar k)=1
\]
for all $k\in K'_\C$.
\end{spherical4}

\it{Proof of Theorem \ref{T:spherical}.} Let $a\in K_\mathbb{C}$
and $x\in G_\mathbb{C}/U$. Write $k'=\kappa^a(k)$, then Lemma
\ref{T:spherical2} implies
 \begin{eqnarray*}
         &&\qquad  E_\q(\phi : ax)\\
         &&=\int_K\tau|_K(k^{-1})(\dqr\phi)(H_\q(kax))\,dk\\
         &&=\int_K\tau|_K((\kappa^{a^{-1}}(k'))^{-1})(\dqr\phi)((I(k'a^{-1}))^{-1}(H_\q(k'x)))\,dk \\
 \end{eqnarray*}
Now use Proposition \ref{T:jacobian}, we have
\begin{eqnarray*}
         &&\qquad  E_\q(\phi : ax)\\
         &&=\int_K\delta(I(ka^{-1}))\tau|_K((\kappa^{a^{-1}}(k))^{-1})
            (\dqr\phi)((I(ka^{-1}))^{-1}H_\q(kx))\,dk.
\end{eqnarray*}
Write $I'=(I(ka^{-1}))^{-1}$ and $k''=(\kappa^{a^{-1}}(k))^{-1}$
for simplicity. By Lemma \ref{T:spherical3}, we have
\begin{eqnarray*}
          &&\quad E_\q(\phi : ax)\\
          &&=\int_K\delta(I'^{-1})\tau|_K(k'')
            \delta_{\n\cap\k}(I')\delta_{\bar \n\cap \k}(\bar I')\tau|_{K_\C}(I')(\dqr\phi)(H_\q(kx))\,dk\\
          &&=\int_K (\delta(I'^{-1})\delta_{\n\cap\k}(I')\delta_{\bar \n\cap \k}(\bar I'))
            \tau|_{K_\C}(k''I')(\dqr\phi)(H_\q(kx))\,dk.
\end{eqnarray*}

Notice that $k''I'=ak^{-1}n$ for some $n\in N_c$, and $N_c$ fix the
values of $\phi$ under the action $\tau|_{K_\C}$. Therefore Lemma
\ref{T:spherical4} implies
\[
   E_\q(\phi : ax)=\int_K \tau|_{K_\C}(ak^{-1})(\dqr\phi)(H_\q(kx))\,dk=\tau|_{K_\C}(a)E_\q(\phi : x).
\]
\qed

\section{Differential equations satisfied by the integral}

\subsection{The differential equations}

From now on throughout this paper, we assume $\q$ is a $\theta$
stable parabolic unless otherwise mentioned. We recall an algebra
homomorphism
\[
  \xi:\U(\g)^K\rightarrow \U(\g')^{K'}
\]
from \cite{V2}.

By PBW theorem,
\[
  \U(\g)=\U(\g')\oplus(\n\U(\g)+\U(\g)\bar \n).
\]
It is known that
\[
  \U(\g)^K\subset\U(\g')^{K'}\oplus(\n\U(\g)\cap\U(\g)\bar \n).
\]
Denote the projection to the first factor by
\[
  \tilde \xi:\U(\g)^K\rightarrow \U(\g')^{K'}.
\]
It is also known that $\tilde \xi$ is an algebra homomorphism. Let
\[
  \eta_\q:\U(\g')\rightarrow \U(\g')
\]
be the algebra homomorphism so that
\[
   \eta_\q(X)=X+\delta_\n(X)
\]
for all $X\in \g'$. $\eta_\q$ maps $\U(\g')^{K'}$ onto itself.
Denote by
\[
  \xi=\eta_\q\circ\tilde \xi:\U(\g)^K\rightarrow \U(\g')^{K'}.
\]

Recall that there is a unique anti-automorphism
\[
  ^{\ve} : \U(\g)\rightarrow \U(\g)
\]
so that
\[
  X^{\ve}=-X, \qquad X\in \g.
\]
Notice that the map "$\phantom{X}^{\ve}$" preserves both
$\U(\g)^K$ and $\U(\g')^{K'}$. If $A$, $B$ are two subalgebras of
$\U(\g)$ which are stable under $^\ve$ and $\eta: A\rightarrow B$
is any algebra homomorphism, denote by $\eta^\ve:A \rightarrow B$
the algebra homomorphism making the following diagram commutes
\[
   \begin{CD}
          A       @>\eta>>              B \\
          @V\ve VV                             @VV\ve V\\
          A       @>\eta^\ve>>          B.
   \end{CD}
\]

Let $\V$, $\tau$, $\V'$ and $\tau'$ be as in section \ref{S:Eint}.
Recall that we have an action $T$ of $\U(\g)\otimes\U(\g)$ on
$\con^\infty(G_\C/U;\V)$ and $\con^\infty(G;\V)$. Denote by $T'$
the analogous action of $\U(\g')\otimes\U(\g')$ on
$\con^\infty(G'_\C/U';\V')$ and $\con^\infty(G';\V')$. It's
obvious that $T_{\U(\g)^K\otimes \U(\g)^K}$ stabilizes
$\con^\infty(G_\C/U;\tau)$, and $T'_{\U(\g')^{K'}\otimes
\U(\g')^{K'}}$ stabilizes $\con^\infty(G'_\C/U';\tau')$. The main
result of this section is:

\begin{dequation}\label{T:dequation}
       For all $X\otimes Y\in \U(\g)^K\otimes\U(\g)^K$ and
       $\phi\in\con^\infty(G'_\C/U';\tau')$,
       \[
          T_{X\otimes Y}E_\q(\phi)=E_\q(T'_{\xi^\ve(X)\otimes \xi(Y)}\phi).
       \]
\end{dequation}

\subsection{Proof of Theorem \ref{T:dequation}}
Notice that we also have
\[
  \U(\g)^K\subset\U(\g')^{K'}\oplus(\bar \n\U(\g)\cap\U(\g)\n).
\]
The projection to the first factor is just the algebra homomorphism
\[
  {\tilde \xi}^\ve:\U(\g)^K\rightarrow \U(\g')^{K'}.
\]

\begin{dequation2}\label{T:dequation2}
Let $X\otimes Y\in\U(\g)^K\otimes\U(\g)^K$ and
$f\in\con^\infty(G_\C/U; \V)$. If $f$ is $N$ invariant, then
\[
  T_{X\otimes Y}f=T_{{\tilde \xi}^\ve(X)\otimes \tilde \xi(Y)}f.
\]
\end{dequation2}

\proof Let $X=X'+X''$, $Y=Y'+Y''$, where $X'={\tilde \xi}^\ve(X)$,
$X''\in U(\g)\n$, $Y'=\tilde \xi(Y)$, $Y''\in U(\g)\bar \n$. Notice
that under the complexification
\[
  \begin{array}{lll}
      1\times\bar{\theta}: &  \g &\rightarrow \g\times\g,\\
                     &  v  &\mapsto (v,\bar{\theta}(v)),
  \end{array}
\]
The Lie algebra of $N$ has a complexification $\n\times \bar \n$.
Therefore the $N$ invariance of $f$ implies $T_{\n\times \bar
\n}f=0$, i.e., $T_{\n\otimes 1+1\otimes \bar \n}f=0$ . So we have
\[
  T_{X\otimes Y}f=T_{X'\otimes Y'}f=T_{{\tilde \xi}^\ve(X)\otimes \tilde \xi(Y)}f.
\]
\qed

\begin{dequation3}\label{T:dequation3}
Let $X\otimes Y\in\U(\g)^{K}\otimes\U(\g)^{K}$ and
$\phi\in\con^\infty(G'_\C/U'; \V')$. Then
\[
  T'_{{\tilde \xi}^\ve(X)\otimes \tilde \xi(Y)}(\dqr\phi)=\dqr T'_{\xi^\ve(X)\otimes
  \xi(Y)}(\phi).
\]
\end{dequation3}

\proof
Let
\[
  \eta_{\bar \q}:\U(\g')\rightarrow \U(\g')
\]
be the algebra homomorphism so that $\eta_{\bar
\q}(X)=X+\delta_{\bar \n}(X)$ for all $X\in \g'$. It's easy to see
that $\eta_{\bar \q}={\eta_\q}^\ve$ and hence
\[
  {\xi}^\ve=(\eta_\q\circ \tilde \xi)^\ve=\eta_{\bar \q}\circ
  {\tilde\xi}^\ve.
\]

Define an algebra homomorphism
\[
  \begin{array}{rl}
         \eta_0:\U(\g'\times\g')=\U(\g')\otimes\U(\g')&\rightarrow
         \U(\g'\times\g')=\U(\g')\otimes\U(\g'),\\
         X'\otimes Y'&\mapsto\eta_{\bar \q}(X')\otimes
         \eta_\q(Y').
  \end{array}
\]
We only need to show that
\[
  T'_Z(\dqr\phi)=\dqr T'_{\eta_0(Z)}(\phi)
\]
for all $Z\in \U(\g'\times\g')$. As $\U(\g'\times\g')$ is generated
by $\g'\times\g'$, it is sufficient to show the above equality holds
for $Z\in \g'\times\g'$. Now assume $Z=(X',Y')\in \g'\times\g'$.
Then $T'_Z$ is an action defined by a vector field on $G'_\C/U'$.
Hence
\[
  T'_Z(\dqr\phi)=T'_Z(\dqr)\phi + \dqr T'_Z(\phi).
\]
We easily find that
\[
  T'_Z(\dqr)=(\delta_{\bar \n}(X')+\delta_\n(Y')) \dqr.
\]\
Therefore
\begin{eqnarray*}
  &&\quad T'_Z(\dqr\phi)\\
  &&=T'_Z(\dqr)\phi + \dqr T'_Z(\phi)\\
  &&=\dqr((\delta_{\bar \n}(X')+\delta_\n(Y'))\phi+T'_Z\phi)\\
  &&=\dqr T'_{\eta_0(Z)}(\phi).
\end{eqnarray*}
\qed

\it{Proof of Theorem \ref{T:dequation}.} Recall that $T_k:
G_\C/U\rightarrow G_\C/U$ is the translation by $k$ for all $k\in
K$. Let $\epsilon_\V$ be the embedding of $\V'$ into $\V$. Write
\[
  f=\epsilon_\V\circ(\dqr\phi)\circ H_\q\in \con^\infty(G_\C/U;\V).
\]
Then
\[
  E_\q(\phi)=\int_K\tau|_K(k^{-1})\circ f\circ T_k\,dk.
\]
Notice that
\[
  T_{X\otimes Y}: \con^\infty(G_\C/U;\V)\rightarrow \con^\infty(G_\C/U;\V)
\]
is a differential operator. Therefore we have
\begin{eqnarray*}
  T_{X\otimes Y}E_\q(\phi)&&=\int_K T_{X\otimes Y}(\tau|_K(k^{-1})\circ f\circ
  T_k)\,dk\\
  &&=\int_K \tau|_K(k^{-1})\circ T_{X\otimes Y}(f\circ T_k)\,dk.
\end{eqnarray*}
$X\otimes Y\in \U(\g)^K\otimes\U(\g)^K$ implies
\[
  T_{X\otimes Y}(f\circ T_k)=T_{X\otimes Y}(f)\circ T_k.
\]
Hence
\[
  T_{X\otimes Y}E_\q(\phi)=\int_K \tau|_K(k^{-1})\circ T_{X\otimes Y}(f)\circ T_k\,dk.
\]
Notice that the map
\[
  \begin{array}{rl}
         \con^\infty(G'_\C/U';\V') &\rightarrow
         \con^\infty(G_\C/U;\V),\\
         \psi &\mapsto \epsilon_\V\circ\psi\circ H_\q
  \end{array}
\]
is $G'_\C$ equivariant. Therefore
\begin{equation}\label{E:equivariant}
  T_{X'\otimes Y'}(\epsilon_\V\circ\psi\circ H_\q)=\epsilon_\V\circ
  T'_{X'\otimes Y'}\psi\circ H_\q
\end{equation}
for all $\quad X'\otimes Y'\in \U(\g')\otimes \U(\g')$,
$\psi\in\con^\infty(G'_\C/U';\V')$.

Now we have
\[
  \begin{array}{ll}
         \medskip
         \quad T_{X\otimes Y}(f)  & \\
         \medskip
         =T_{{\tilde \xi}^\ve(X)\otimes \tilde \xi(Y)}f & \qquad\textrm{Lemma \ref{T:dequation2}} \\
         \medskip
         =T_{{\tilde \xi}^\ve(X)\otimes \tilde \xi(Y)}(\epsilon_\V\circ(\dqr\phi)\circ
         H_\q) &\\
         \medskip
         =\epsilon_\V\circ (T'_{{\tilde \xi}^\ve(X)\otimes \tilde \xi(Y)}(\dqr\phi))\circ
         H_\q &  \qquad\textrm{By \eqref{E:equivariant}}\\
         =\epsilon_\V\circ (\dqr T'_{\xi^\ve(X)\otimes \xi(Y)}\phi)\circ
         H_\q. & \qquad\textrm{Lemma \ref{T:dequation3}}
  \end{array}
\]
In conclusion,
\begin{eqnarray*}
        &&\quad T_{X\otimes Y}E_\q(\phi)\\
        &&=\int_K \tau|_K(k^{-1})\circ T_{X\otimes Y}(f)\circ T_k\,dk\\
        &&=\int_K \tau|_K(k^{-1})\circ \epsilon_\V \circ (\dqr T'_{\xi^\ve(X)\otimes \xi(Y)}\phi)\circ
           H_\q\circ T_k\,dk\\
        &&=E_\q(T'_{\xi^\ve(X)\otimes \xi(Y)}\phi).
\end{eqnarray*}
\qed

\section{Actions of $\U(\g)^K$ on bottom layers of cohomologically induced
modules}

\subsection{Zuckerman's functors}

We call a vector space a weak $(\U(\g)^K,K)$ module if it is a
locally finite $K$ module as well as a $\U(\g)^K$ module so that the
$K$ action and the $\U(\g)^K$ action commute.  Notice that every
$(\g,K)$ module is a weak $(\U(\g)^K,K)$ module. In this section, we
study the weak $(\U(\g)^K,K)$ module structure of a bottom layer of
a cohomologically induced module. We recall some notations from
\cite{KV} for cohomological inductions and bottom layer maps.

$\Gamma^0$ is the Zuckerman functor from the category of $(\g,K')$
modules to the category of $(\g,K)$ modules. It's given by
\begin{eqnarray*}
         &&\Gamma^0(M)=\textrm{\, sum of all finite dimensional $\k$
         invariant subspaces}\\
         &&\qquad \qquad \textrm{of $M$ for which the action of $\k$
          globalize to $K$\,}.
\end{eqnarray*}
This is a left exact covariant functor. Let $\Gamma^i$ be its
$i^{th}$ derived functor, $i=0,1,2,\cdots$. $\Gamma^0_K$ is the
Zuckerman functor from the category of $(\k,K')$ modules to the
category of $(\k,K)$ modules. $\Gamma^i_K$, $i=0,1,2,\cdots$, are
its derived functors.

$\pro_{\q,K'}^{\g,K'}$ is the exact covariant functor from the
category of $(\q,K')$ modules to the category of $(\g,K')$ modules
given by
\[
   \pro_{\q,K'}^{\g,K'}(M)=\Hom_{\U(\q)}(\U(\g),M)_{K'},
\]
where the $\U(\g)$ action on $\Hom_{\U(\q)}(\U(\g),M)$ is given by
\[
  (X\phi)(u)=\phi(uX),\qquad X,u\in \U(\g),\,\phi\in
  \Hom_{\U(\q)}(\U(\g),M),
\]
the $K'$ action on $\Hom_{\U(\q)}(\U(\g),M)$ is given by
\[
  (k\phi)(u)=k(\phi(\Ad_{k^{-1}}u)),\qquad k\in K',\, u\in
  \U(\g),\,\phi\in \Hom_{\U(\q)}(\U(\g),M),
\]
and $\Hom_{\U(\q)}(\U(\g),M)_{K'}$ means the $K'$ finite vectors of
$\Hom_{\U(\g')}(\U(\g),M)$. $\pro_{\q\cap\k,K'}^{\k,K'}$ is the
analogous exact covariant functor from the category of
$(\q\cap\k,K')$ modules to the category of $(\k,K')$ modules. $\FF$
is used to denote forgetful functors. For example,
$\FF_{\g',K'}^{\q,K'}$ is the forgetful functor from the category of
$(\g',K')$ modules to the category of $(\q,K')$ modules via trivial
$\n$ action.

Let $M'$ be a $(\g',K')$ module and define another $(\g',K')$ module
\[
  M'^\#=M'\otimes \wedge^{\operatorname{top}}\n.
\]
Define
\[
  \RR^i(M')=(\Gamma^i \circ \pro_{\q,K'}^{\g,K'} \circ
  \FF_{\g',K'}^{\q,K'})(M'^\#).
\]
$\RR^i$ is the cohomological induction functor from the category of
$(\g',K')$ modules to the category of $(\g,K)$ modules. Define
\[
  \RR_\circ^i(M')=(\Gamma_K^i \circ\FF_{\g,K'}^{\k,K'}
  \circ \pro_{\q,K'}^{\g,K'} \circ\FF_{\g',K'}^{\q,K'})(M'^\#).
\]
Notice that  $\U(\g)^K$ acts on $(\FF_{\g,K'}^{\k,K'}\circ
\pro_{\q,K'}^{\g,K'} \circ \FF_{\g',K'}^{\q,K'})(M'^\#)$ by
$(\k,K')$ module endomorphisms. Therefore by using the functor
$\Gamma_K^i$, we have an action of $\U(\g)^K$ on $\RR_\circ^i(M')$
by $(\k,K)$ module endomorphisms. Define
\[
  \RR_K^i(M')=(\Gamma_K^i \circ\pro_{\q\cap\k,K'}^{\k,K'} \circ
  \FF_{\g',K'}^{\q\cap\k,K'})(M'^\#).
\]
Notice that $\U(\g')^{K'}$ acts on
$\FF_{\g',K'}^{\q\cap\k,K'}(M'^\#)$ by $(\q\cap\k,K')$ module
endomorphisms. Therefore by applying the functor
$\pro_{\q\cap\k,K'}^{\k,K'}$, we have an action of $\U(\g')^{K'}$ on
$(\pro_{\q\cap\k,K'}^{\k,K'} \circ
\FF_{\g',K'}^{\q\cap\k,K'})(M'^\#)$ by $(\k,K')$ module
endomorphisms. Then by applying $\Gamma_K^i$, we have an action of
$\U(\g')^{K'}$ on $\RR_K^i(M')$ by $(\k,K)$ module endomorphisms.

\begin{bottom}\label{T:bottom}
For any $(\g',K')$ module $M'$, $\RR^i(M')$ is canonically
isomorphic to $\RR_\circ^i(M')$ as weak $(\U(\g)^K,K)$ modules.
\end{bottom}

\proof
Write $M''=(\pro_{\q,K'}^{\g,K'} \circ
\FF_{\g',K'}^{\q,K'})(M'^\#)$ and let
\[
   0\rightarrow M'' \rightarrow J_0 \rightarrow J_1 \rightarrow J_2
   \rightarrow \cdots
\]
be an injective resolution of $M''$ in the category of $(\g,K')$
modules. By Proposition 2.57 of \cite{KV}, the exact functor
$\FF_{\g,K'}^{\k,K'}$ sends injectives to injectives. Therefore
\[
   0\rightarrow \FF_{\g,K'}^{\k,K'}(M'') \rightarrow
   \FF_{\g,K'}^{\k,K'}(J_0) \rightarrow
   \FF_{\g,K'}^{\k,K'}(J_1) \rightarrow
   \FF_{\g,K'}^{\k,K'}(J_2) \rightarrow \cdots
\]
is an injective resolution of $\FF_{\g,K'}^{\k,K'}(M'')$ in the
category of $(\k,K')$ modules.

Notice that
\[
  0\rightarrow \Gamma^0(J_0) \rightarrow \Gamma^0(J_1)
  \rightarrow \Gamma^0(J_2) \rightarrow \cdots
\]
and
\[
  0\rightarrow \Gamma_K^0(\FF_{\g,K'}^{\k,K'}(J_0))
  \rightarrow \Gamma_K^0(\FF_{\g,K'}^{\k,K'}(J_1)) \rightarrow
  \Gamma_K^0(\FF_{\g,K'}^{\k,K'}(J_2)) \rightarrow \cdots
\]
are exactly the same as sequences of weak $(\U(\g)^K,K)$ modules.
Take the $i^{th}$ cohomology of both sequences, we get
$\RR^i(M')=\RR_\circ^i(M')$ as weak $(\U(\g)^K,K)$ modules.
\qed

\subsection{Bottom layer maps}

Let $M'$ be a $(\g',K')$ module. Define a $(\k,K')$ module
homomorphism
\[
   \beta_{M'} :(\FF_{\g,K'}^{\k,K'}\circ \pro_{\q,K'}^{\g,K'} \circ
   \FF_{\g',K'}^{\q,K'})(M'^\#) \rightarrow
   (\pro_{\q\cap\k,K'}^{\k,K'} \circ
   \FF_{\g',K'}^{\q\cap\k,K'})(M'^\#)
\]
by
\[
  (\beta_{M'}(\phi))(r)=\phi(r),\qquad \phi\in
  (\FF_{\g,K'}^{\k,K'}\circ \pro_{\q,K'}^{\g,K'} \circ
  \FF_{\g',K'}^{\q,K'})(M'^\#),\, r\in \U(\k).
\]

\begin{bottom2}\label{T:bottom2}
For all $X\in \U(\g)^K$, the following diagram commutes
\[
  \begin{CD}
        (\FF_{\g,K'}^{\k,K'}\circ \pro_{\q,K'}^{\g,K'} \circ
        \FF_{\g',K'}^{\q,K'})(M'^\#) @>\beta_{M'}>>
        (\pro_{\q\cap\k,K'}^{\k,K'} \circ
        \FF_{\g',K'}^{\q\cap\k,K'})(M'^\#)\\
        @V{X}VV        @V\tilde \xi (X) VV                \\
        (\FF_{\g,K'}^{\k,K'}\circ \pro_{\q,K'}^{\g,K'} \circ
        \FF_{\g',K'}^{\q,K'})(M'^\#) @>\beta_{M'}>>
        (\pro_{\q\cap\k,K'}^{\k,K'} \circ \FF_{\g',K'}^{\q\cap\k,K'})(M'^\#).\\
  \end{CD}
\]
\end{bottom2}

\proof
Let
\[
   X=aX_1+X_2,
\]
where $a\in \n$, $X_1\in \U(\g)$, and
\[
   X_2=\tilde \xi (X)\in \U(\g')^{K'}.
\]
Let
\[
   \phi\in (\FF_{\g,K'}^{\k,K'}\circ \pro_{\q,K'}^{\g,K'}
   \circ \FF_{\g',K'}^{\q,K'})(M'^\#)\subset
   \Hom_{\U(\q)}(\U(\g),M'^\#),
\]
and $r\in \U(\k)$. We have
\begin{eqnarray*}
         &&\quad(\beta_{M'}(X\phi))(r)\\
         &&=(X\phi)(r)=\phi(rX)=\phi(Xr)\\
         &&=\phi(aX_1 r)+\phi(X_2 r)\\
         &&=a(\phi(X_1))+X_2(\phi(r))\\
         &&=X_2(\beta_{M'}(\phi)(r))\\
         &&=\tilde \xi(X)(\beta_{M'}(\phi))(r).
\end{eqnarray*}
Therefore
\[
  \beta_{M'}(X\phi)=\tilde \xi(X)(\beta_{M'}(\phi)).
\]
\qed

Apply the functor $\Gamma_K^i$ to $\beta_{M'}$, we define the bottom
layer map
\[
  \BB_{M'}=\Gamma_K^i(\beta_{M'}) : \RR_\circ^i(M')\rightarrow
  \RR_K^i(M').
\]

Apply $\Gamma_K^i$ to the commutative diagram of the above Lemma, we
get

\begin{bottom3}\label{bottom3}
For all $X\in \U(\g)^K$, the following diagram commutes:
\[
  \begin{CD}
         \RR_\circ^i(M') @>\BB_{M'}>> \RR_K^i(M')\\
         @VXVV           @V \tilde \xi(X) VV                \\
         \RR_\circ^i(M') @>\BB_{M'}>> \RR_K^i(M').\\
  \end{CD}
\]
\end{bottom3}

\subsection{Bottom layers}\label{S:Blayers}

Let $\sigma$ be a finite dimensional continuous representation of
$K$ on $W$. Let
\begin{equation}\label{E:Wprime}
      W'=\{\,v\in W \mid \sigma(\n\cap\k)v=0\,\}.
\end{equation}
Define a representation $\sigma'$ of $K'$ on $W'$ by
\begin{equation}\label{E:sigmaprime}
     \sigma'(k)v=\delta_{\bar \n\cap\p}(k)\sigma(k)v, \qquad k\in K',\,
     v\in W.
\end{equation}
We call $\sigma'$ the representation $\q$-associated to $\sigma$. If
$\sigma$ is irreducible of type $\alpha$, then $\sigma'$ is also
irreducible. We call the equivalent class of $\sigma'$ the $K'$ type
$\q$-associated to $\alpha$.

From now on in this section, we fix a $K$ type $\alpha$. We also fix
an irreducible representation $\sigma_0$ of $K$ of type $\alpha$ on
a vector space $W_\alpha$. Let $\sigma_0'$ be the representation of
$K'$ on $W_{\alpha'}$ which is $\q$-associated to $\sigma_0$, and
$\alpha'$ be the $K'$ type $\q$-associated to $\alpha$.

Denote by $S=\dim(\n\cap\k)$ as usual. The most interesting case of
cohomological induction is when $i=S$. Denote by $M=\RR^S(M')$ and
$L=\Hom_{K'}(W_{\alpha'},M')$. $L$ is a $\U(\g')^{K'}$ module by the
action on $M'$. We make $L$ a $\U(\g)^{K}$ module by the formula
\[
  Xv=\xi(X)v, \qquad X\in \U(\g)^K, v\in L.
\]
$L\otimes W_\alpha$ is a weak $(\U(\g)^K,K)$ module by the action of
$\U(\g)^K$ on the first factor, and the action of $K$ on the second
factor.

\begin{bottom4}
Let $M'$ be a $(\g',K')$ module, $M=\RR^S(M')$. $\alpha$ is said to
be in the bottom layer of $M$ if the map induced by $\BB_{M'}$,
\[
  \BB_{M'}(\alpha) : M(\alpha)=\RR_\circ^S(M')(\alpha)\rightarrow
  \RR_K^S(M')(\alpha)
\]
is bijective and nonzero.
\end{bottom4}

We know that under a weak condition on $M'$, $\alpha$ is in the
bottom layer of $M$ if and only if $M'(\alpha')\neq 0$. The
condition holds when $M'$ has an infinitesimal character (\cite{KV},
Theorem 5.80, Corollary 5.72). I don't know whether the following
result is new or not.

\begin{bottom5}\label{T:bottom5}
With the notations as above. If $\alpha$ is in the bottom layer of
$M$, then $M(\alpha)$ is isomorphic to $L\otimes W_\alpha$ as a weak
$(\U(\g)^K,K)$ module.
\end{bottom5}

\proof
We easily find that the following diagram commutes
\[
  \begin{CD}
         M'^\#  @=        M'\otimes  \wedge^{\operatorname{top}}\n\\
         @VXVV           @V \eta_\q(X)\otimes 1 VV                \\
         M'^\# @=        M'\otimes  \wedge^{\operatorname{top}}\n\\
  \end{CD}
\]
for all $X\in \U(\g')$. This implies another commutative diagram
\[
  \begin{CD}
       \Hom_{K'}(W_{\alpha'}\otimes \wedge^{\operatorname{top}}\n,M'^\#)   @=  \Hom_{K'}(W_{\alpha'},M') \\
         @VXVV           @V \eta_\q(X)VV                \\
       \Hom_{K'}(W_{\alpha'}\otimes \wedge^{\operatorname{top}}\n,M'^\#)   @=  \Hom_{K'}(W_{\alpha'},M') \\
  \end{CD}
\]
for all $X\in \U(\g')^{K'}$. Since
\[
   \operatorname{H}^0(\n\cap\k, W_\alpha)\otimes
       \wedge^{\operatorname{top}}(\n\cap\k) \cong
   W_{\alpha'}\otimes \wedge^{\operatorname{top}}\n
\]
as $K'$ modules, Theorem 4.155 of \cite{KV} gives an isomorphism
\[
  \Hom_K(W_\alpha, \RR_K^S(M'))\cong \Hom_{K'}(
  W_{\alpha'}\otimes \wedge^{\operatorname{top}}\n,M'^\#).
\]
The naturalness implies this is an isomorphism of $\U(\g')^{K'}$
modules. Therefore Proposition \ref{bottom3} and the above
commutative diagram implies we have a commutative diagram
\[
  \begin{CD}
         \Hom_K(W_\alpha, \RR^S(M'))  @>>> L=\Hom_{K'}(W_{\alpha'},M') \\
         @VXVV           @V  \xi(X) VV                \\
         \Hom_K(W_\alpha, \RR^S(M'))  @>>> L=\Hom_{K'}(W_{\alpha'},M') \\
  \end{CD}
\]
for all $X\in \U(\g)^K$. The horizontal maps are bijective as
$\alpha$ is in the bottom layer. This finishes the proof.
\qed

\section{Matrix coefficient of a bottom layer}\label{S:matrix}

\subsection{The matrix coefficient of a given $K$ type}

Let $W$ be a finite dimensional weak $(\U(\g)^K,K)$ module. Denote
by $\sigma$ the actions of $\U(\g)^K$ and $K$ on $W$. Define
$\aut_W(G)$ to be the space of all functions
$\phi\in\con^\infty(G;\End_\C(W))$ so that
\[
  \begin{array}{l}
        \phi(k_1 x k_2)=\sigma(k_1)\circ \phi(x)\circ \sigma(k_2),\qquad
        k_1,k_2\in K,\,x\in G, \textrm{  and }\\
        T_{X\otimes Y}\phi(x)=\sigma(X^\ve)\circ \phi(x) \circ
        \sigma(Y), \qquad x\in G,\,X\otimes Y\in \U(\g)^K\otimes \U(\g)^K.
 \end{array}
\]
Notice that every function in $\aut_W(G)$ is real analytic.

If $M$ is an admissible $(\g,K)$ module, denote by $M^{\oc}$ the
contragredient $(\g,K)$ module of $M$. If $V$ is a finite
dimensional continuous $K$ module and $\alpha$ is a $K$ type,
define $V^{\oc}$ and $\alpha^{\oc}$ similarly. Let's recall the
basic fact of the existence and uniqueness of the matrix
coefficient map.

\begin{eandu}(\cite{WD}, Theorem 8.7)\label{T:eandu}
If $M$ is an admissible $(\g,K)$ module, then there is a unique
$(\g\times\g,K\times K)$ intertwining linear map
\[
  c_M : M^{\oc}\otimes M\rightarrow \con^\infty(G;\C)
\]
such that
\[
   c_M(u\otimes v)(1)=\langle u,v\rangle, \qquad u\otimes v
   \in M^c \otimes M,
\]
where $1$ is the identity element of $G$.
\end{eandu}

Now let $M$ be an arbitrary finitely generated admissible $(\g,K)$
module, $\alpha$ a $K$ type, and $W=M(\alpha)$. Then $W$ is a weak
$(\U(\g)^K,K)$ module. $W^{\oc} = M^{\oc}(\alpha^{\oc})$ is the
contragredient of the $K$ module $W$. We use $c_{M,\alpha}$ to
denote the function on $G$ with values in $\End_{\C}(W)$ defined
by
\[
  \la u,c_{M,\alpha}(x)(v)\ra=c_M(u\otimes v)(x), \qquad u\in
  W^{\oc}, \,v\in W,\, x\in G,
\]
where $c_M$ is the matrix coefficient map in the above proposition.
The function $c_{M,\alpha}$ can also be described as the following.
Take a Hilbert representation
\[
  \pi_G: G\times M_G\rightarrow M_G
\]
whose underlying $(\g,K)$ module is $M$. Let
\[
  s_W: M_G\rightarrow W
\]
be the continuous linear map so that its restriction to $M_G(\beta)$
is zero for all $K$ type $\beta$ other than $\alpha$, and its
restriction to $W$ is the identity map. Then we have
\[
  c_{M,\alpha}(x)=s_W \circ \pi_G(x)\circ j_W \in \End_{\C}(W), \qquad x\in G,
\]
where $j_W$ is the inclusion of $W$ into $M_G$. By the above
formula, it's not difficult to see that $c_{M,\alpha}\in \aut_W(G)$.
It's obvious that $c_{M,\alpha}(1)$ is the identity map on $W$. The
following proposition says that these two conditions determines
$c_{M,\alpha}$ uniquely.

\begin{coef1}\label{T:coef1}
Let $W$ be a finite dimensional weak $(\U(\g)^K,K)$ module. If $W$
is primary as a $K$ module, then the map $\phi\mapsto \phi(1)$
gives an injection from $\aut_W(G)$ into $\End_{\U(\g)^K,K}(W)$.
\end{coef1}

\proof
Let $\phi\in \mathcal{A}_W(G)$. It's easy to see that
$\phi(1)\in \End_{\U(\g)^K,K}(W)$. Now we assume $\phi(1)=0$.

Assume $W=L\otimes W_0$, where $L$ is a $\U(\g)^K$ module and
$W_0$ is an irreducible $K$ module. Define a map
\[
  \begin{array}{lrl}
      \Tr_0 : &  \End_\C(W)=\End_\C(L)\otimes\End_\C(W_0) &\rightarrow \End_\C(L),\\
                     &  x\otimes y  &\mapsto \Tr(y)x.
  \end{array}
\]
Let $\tilde{\phi}=\Tr_0 \circ \phi$ be an $\End_\C(L)$ valued real
analytic function, and denote by $\sigma_0$ the action of
$\U(\g)^K$ on $L$. Then we have
\begin{equation*}
        \left\{
          \begin{aligned}
            &\tilde{\phi}(1)=0;\\
            &\tilde{\phi}(kxk^{-1})=\tilde{\phi}(x),\qquad \qquad x\in
                G,\,k\in K;\\
            &T_{X\otimes Y}\tilde{\phi}(x)=\sigma_0(X^\ve) \circ
              \tilde{\phi}(x)\circ \sigma_0(Y),\qquad x\in
              G,\,X\otimes Y\in \U(\g)^K\otimes \U(\g)^K.
         \end{aligned}
      \right.
\end{equation*}

We essentially copy the following arguments from \cite{FJ} to
deduce that $\tilde{\phi}=0$. Let $X\in \U(\g)$. Denote by
\[
  \tilde{X}=\int_K \Ad_k X \,dk\,\in \U(\g)^K.
\]
We have
\[
  T_{X\otimes 1}\tilde{\phi}(1)=\Ad_k(T_{X\otimes 1}\tilde{\phi})(1)
  =T_{\Ad_k(X)\otimes 1}(\Ad_k\tilde{\phi})(1)=T_{\Ad_k(X)\otimes 1}\tilde{\phi}(1)
\]
for all $k\in K$. Integrate the above equality over $K$ we get
\[
  T_{X\otimes 1}\tilde{\phi}(1)=T_{\tilde{X}\otimes 1}\tilde{\phi}(1)
  =\sigma(\tilde{X}^\ve)\circ (\tilde{\phi}(1))=0.
\]
Therefore $\tilde{\phi}=0$ since it is real analytic and the group
$G$ is connected.

Let $\V_0=\Span\{\,\phi(x)\mid x\in G\,\}$ be a subspace of
\[
   \End_\C(W)=\End_\C(L)\otimes\End_\C(W_0).
\]
Since $\V_0$ is $K\times K$ stable,
\[
  \V_0=\V_1\otimes \End_\C(W_0)
\]
for some subspace $\V_1$ of $\End_\C(L)$. Now it's clear that
$\tilde{\phi}=0$ implies $\V_1=0$, and consequently $\phi=0$.
\qed

\subsection{Matrix coefficient of a bottom layer}\label{S:Mlayer}
Now we are able to get the result on the matrix coefficients of a
bottom layer of a cohomologically induced module. Let $M=\RR^S(M')$
be a cohomologically induced module, where $S=\dim(\n\cap\k)$ and
$M'$ is a finitely generated admissible $(\g',K')$ module. Let
$\alpha$, $\alpha'$, $W_\alpha$, $W_{\alpha'}$ and $L$ be as in
section \ref{S:Blayers}. Assume $\alpha$ is in the bottom layer of
$M$. Denote by
\[
    W'=M'(\alpha')=L\otimes W_{\alpha'}
\]
as a weak $(\U(\g')^{K'},K')$ module. Denote by
\[
   W=L\otimes W_\alpha
\]
as a weak $(\U(\g)^K,K)$ module. Proposition \ref{T:bottom5} says
that $W\cong M(\alpha)$ as weak $(\U(\g)^{K},K)$ modules. We fix
an isomorphism and identify $W$ with $M(\alpha)$. Denote by
$\sigma$ the actions of $\U(\g)^K$ and $K$ on $W$; and $\sigma'$
the actions of $\U(\g')^{K'}$ and $K'$ on $W'$.

Define a representation $\tau$ of $K\times K$ on $\End_\C(W)$ by
\[
  \tau(k,l)(f)=\sigma(k)\circ f\circ \sigma(l^{-1}),\qquad k,l\in K,\, f\in
  \End_\C(W),
\]
and a representation $\tau'$ of $K'\times K'$ on $\End_\C(W')$ by
\[
  \tau'(k,l)(f)=\sigma'(k)\circ f\circ \sigma'(l^{-1}),\qquad k,l\in K',\, f\in
  \End_\C(W').
\]
We have a decomposition
\[
   W=W'\oplus\sigma(\bar{\n}\cap \k)W.
\]
Therefore we may view every element of $\End_\C(W')$ as an element
of $\End_\C(W)$ which vanishes on $\sigma(\bar{\n}\cap \k)W$. Denote
by
\[
  j_\q:\End_\C(W')\rightarrow \End_\C(W)
\]
the corresponding embedding.

\begin{qassociated}
$j_\q$ induces an isomorphism from $\tau'$ to the representation
$\q$-associated to $\tau$.
\end{qassociated}

We omit the straight forward proof.

We have defined the matrix coefficient $c_{M,\alpha}\in\aut_W(G)$.
Similarly we have the matrix coefficient $c_{M',\alpha'}\in
\aut_{W'}(G')$. Proposition \ref{T:coef1} tells us that
$c_{M,\alpha}$ is the unique element in $\aut_W(G)$ such that
\[
  c_{M,\alpha}(1)=1_W.
\]
The same holds for $c_{M',\alpha'}$. Notice that
\[
  \aut_W(G)\subset \con^\infty(G;\tau) \quad \textrm{and}\quad
  \aut_{W'}(G')\subset \con^\infty(G';\tau').
\]
Denote by $d_{M,\alpha}\in\con^\infty(G_\C/U;\tau)$ the function
corresponds to $c_{M,\alpha}$ in Proposition \ref{T:fj2}.
$d_{M',\alpha'}\in\con^\infty(G'_\C/U';\tau')$ is defined
similarly.

The main theorem of this section is the following:

\begin{coef4}\label{T:coef4}
With the notations as above, we have
\begin{equation}\label{E:coef4}
  d_{M,\alpha}=\frac{\degree(\alpha)}{\degree(\alpha')}E_\q(j_\q\circ d_{M',\alpha'}).
\end{equation}
\end{coef4}

\subsection{A proof of Theorem \ref{T:coef4}}

We need the following elementary lemma.

\begin{coef2}\label{T:coef2}
Let $\sigma$ be a continuous representation of $K$ on a finite
dimensional vector space $W$. Let $P\in \End_\C(W)$ be a linear map
which stabilize every irreducible submodule of $W$. Suppose that
$P^2=P$ and there is a constant $c$ such that
\[
  \frac{\dim(P(W_0))}{\dim(W_0)}=c
\]
for every irreducible submodule $W_0$ of $W$. Then
\[
  \int_K \sigma(k^{-1})\circ P\circ \sigma(k)\, dk=c1_W.
\]
\end{coef2}

\proof We assume $W$ is irreducible without lose of generality.
Denote by $C$ the integral in the lemma. It's clear that $C$ is an
intertwining operator. The irreducibility of $W$ implies that
$C=c'1_W$ for some constant $c'$. By taking trace of both side of
the above equality we get $c'=c$.
\qed

Now we are going to prove the theorem. It's clear from the
definition of $E_\q$ that
\[
  E_\q(j_\q\circ d_{M',\alpha'})(1U)=\int_K \sigma(k^{-1}) \circ j_\q(1_{W'})\circ
  \sigma(k) \,dk
\]
$j_\q(1_{W'})$ satisfies the conditions on $P$ in Lemma
\ref{T:coef2}, with the constant
$c=\frac{\degree(\alpha')}{\degree(\alpha)}$. Therefore Lemma
\ref{T:coef2} implies that
\[
  E_\q(j_\q\circ d_{M',\alpha'})(1U)=c1_W.
\]
This proves that both hand side of \eqref{E:coef4} has value $1_W$
at the point $1U$. Denote by
\[
   \psi_0=E_\q(j_\q\circ d_{M',\alpha'}),
\]
and $\phi_0$ the corresponding function in $\con^\infty(G;\tau)$. By
Proposition \ref{T:coef1}, we only need to show that $\phi_0\in
\aut_W(G)$ to finish the proof.

Let $X\otimes Y\in \U(\g)^K\otimes \U(\g)^K$. We easily check that
the following diagram commutes:
\[
   \begin{CD}
          \End_\C(W')          @>j_\q>>         \End_\C(W)  \\
          @V A'VV                               @VV A V\\
          \End_\C(W')          @>j_\q>>         \End_\C(W),
   \end{CD}
\]
where $A'$ is defined by
\[
  f'\mapsto \sigma'(\xi(X^\ve))\circ f'\circ \sigma'(\xi(Y)),
\]
and $A$ is defined by
\[
  f\mapsto \sigma(X^\ve)\circ f\circ \sigma(Y).
\]

As $(\xi^\ve(X))^\ve=\xi(X^\ve)$ and $c_{M',\alpha'}\in
\aut_{W'}(G')$, we have
\[
  T'_{\xi^\ve(X)\otimes \xi(Y)}(c_{M',\alpha'})=A'\circ c_{M',\alpha'}.
\]
Now Lemma \ref{T:doperator} and Lemma \ref{T:fj} implies that
\begin{equation}\label{E:psiprime}
  T'_{\xi^\ve(X)\otimes \xi(Y)}(d_{M',\alpha'})=A'\circ d_{M',\alpha'}.
\end{equation}

We have
\[
   \begin{array}{ll}
      \medskip
      \quad T_{X\otimes Y}(\psi_0) & \\
      \medskip
      =T_{X\otimes Y}E_\q(j_q\circ d_{M',\alpha'}) &\\
      \medskip
      =E_\q(T'_{\xi^\ve(X)\otimes \xi(Y)}(j_\q\circ d_{M',\alpha'})) &
         \qquad \textrm{Theorem } \ref{T:dequation} \\
      \medskip
      =E_\q(j_\q\circ (T'_{\xi^\ve(X)\otimes \xi(Y)}(d_{M',\alpha'}))) & \\
      \medskip
      =E_\q(j_\q \circ A'\circ d_{M',\alpha'}) & \qquad \textrm{By }\eqref{E:psiprime}\\
      \medskip
      =E_\q(A \circ j_\q\circ d_{M',\alpha'}) & \qquad \textrm{By the above commutative diagram}\\
      \medskip
      =A \circ E_\q(j_\q \circ d_{M',\alpha'}) & \qquad A \textrm{ commutes the action of } K\times K\\
      =A\circ \psi_0. &
   \end{array}
\]
Again Lemma \ref{T:doperator} and Lemma \ref{T:fj} implies
\[
  T_{X\otimes Y}\phi_0=A\circ \phi_0.
\]
Therefore $\phi_0\in \aut_W(G)$. This finishes the proof of Theorem
\ref{T:coef4}.

\subsection{An example on $A_\q(\lambda)$
modules}\label{S:aqlambda}

We discuss the case when $\dim(M')=1$ as an example. Let
\[
  \lambda : G'\rightarrow \C^\times
\]
be a continuous character, unitary or not. Still denote by $\lambda$
its holomorphic extension to $G'_\C$ and its differential. Fix a
maximal torus $T_{max}$ of $K'$ with complexified Lie algebra $\t$.
Then $T_{max}$ is also a maximal torus of $K$. Let
\[
  \Lambda=\lambda|_\t+\delta_{\n\cap\p}|_\t.
\]
Let $\alpha$ be the $K$ type with an extremal weight $\Lambda$.
Assume that $\Lambda$ is dominant with respect to $\n\cap \k$. Then
$\lambda|_{K'}$ is the $K'$ type $\q$-associated to $\alpha$, and
$\alpha$ is the unique bottom layer of $A_\q(\lambda)$. Recall that
$A_\q(\lambda)$ is just the cohomologically induced module
$\RR^S(\lambda)$. For every $x\in G$, denote by
\[
 \phi_{\q,\lambda}(x)=\frac{1}{\deg(\alpha)}\Tr(c_{A_\q(\lambda),\alpha}(x)).
\]

Let $\sigma_0$ be an irreducible unitary representation of $K$ on
$W_\alpha$, of type $\alpha$. Fix a unit vector $v_0$ in $W_\alpha$,
of weight $\Lambda$. Let $P_\Lambda\in \End_\C(W_\alpha)$ be the
orthogonal projection onto $\C v_0$. Define a map
\[
  \begin{array}{rl}
                \tilde{H}_\q : G & \rightarrow  G'_\C,\\
                     x   & \mapsto y\overline{\theta(y^{-1})},
  \end{array}
\]
where $y\in G'_\C$ is an element so that $H_\q(xU)=yU'$. Fix an
identification of $A_\q(\lambda)(\alpha)$ with $W_\alpha$ as $K$
modules. Now Theorem \ref{T:coef4} easily implies
\[
    c_{A_\q(\lambda),\alpha}(kx^2)= \deg(\alpha) \int_K (\delta_\n \lambda)(\tilde{H}_\q(lx))
         \sigma_0(kl^{-1})\circ P_\Lambda \circ \sigma_0(l) \,dl\\
\]
for all  $k\in K$, $x\in \exp(\p_0)$. Consequently,
\[
 \phi_{\q,\lambda}(kx^2)=\int_K \la\sigma_0(lkl^{-1})v_0,\, v_0 \ra
      (\delta_\n \lambda)(\tilde{H}_\q(lx))\,dl, \qquad k\in K, x\in \exp(\p_0).
\]
This formula confirms the assumption in Theorem 4.3 of \cite{Li2}
as a special case.

\section{Li positivity}

In this sections, we discuss a positivity property on
representations as an example of the application of our Eisenstein
integral. This positivity property is very useful in the study of
branching rules, for example, in the study of theta
correspondence. All the omitted proofs in this section can be
found in the author's HKUST thesis.

\subsection{The definition of Li positivity}

Let $\sigma$ be a finite dimensional continuous unitary
representation of $K$ on $W$. Define a representation $\tau$ of
$K\times K$ on $\End_\C(W)$ by
\[
  \tau(k,l)(f)=\sigma(k)\circ f\circ \sigma(l^{-1}),\qquad k,l\in K,\, f\in
  \End_\C(W).
\]
Recall that an operator $A\in \End_\C(W)$ is said to be positive
definite if
\[
  \la Av,v\ra\geq 0,\qquad v\in W.
\]

\begin{lipositive}
Let $\phi\in \con^\infty(G;\tau)$. $\phi$ is said to be Li positive
if $\phi(x)$ is positive definite for all $x\in\exp(\p_0)$. Let
$\psi\in \con^\infty(G_\C/U;\tau)$. $\psi$ is said to be Li positive
if $\psi(x)$ is positive definite for all $x\in G_\C/U$.
\end{lipositive}

We use this terminology since such kinds of properties are firstly
used by Jian-Shu Li in the study of theta correspondence \cite{Li1},
\cite{Li2}, \cite{HL}. If $\phi\in \con^\infty(G;\tau)$ is Li
positive, then $\Tr\circ \phi$ satisfies Flensted-Jensen's formula
in Jian-Shu Li's terminology. In some sense, the concept of Li
positivity is dual to that of positive definiteness of functions.

The following lemma will be used later.
\begin{liflensted}\label{T:liflensted}
Let $\phi\in \con^\infty(G;\tau)$ and $\psi\in
\con^\infty(G_\C/U;\tau)$ be a pair of corresponding functions in
Proposition \ref{T:fj2}. Then $\phi$ is Li positive if and only if
$\psi$ is Li positive.
\end{liflensted}

\proof
The "if" part is clear. Assume $\phi$ is Li positive. Let
$x\in G_\C/U$. Recall that we have a complexification square $\S$:
\[
   \begin{CD}
      G/K         @>p_K>>     G  \\
      @Vv_GVV                 @VVu_GV\\
      G_\C/U      @>p_U>>     G_\C.
   \end{CD}
\]
Since
\[
  G_\C/U=K_\C.v_G(G/K),
\]
we have
\[
  x=k.v_G(y) \quad\textrm{for some}\quad k\in K_\C,\,y\in G/K.
\]
Then
\[
  \psi(x)=\tau(k,\bar k)\psi(v_G(y))=\sigma(k)\circ
  \phi(p_K(y))\circ \sigma({\bar k}^{-1}).
\]
Notice that $p_K(y)\in \exp(\p_0)$ and $\sigma({\bar k}^{-1})$ is
the adjoint operator of $\sigma(k)$. Therefore $\psi(x)$ is
positive definite.
\qed

Let $M$ be a finitely generated admissible $(\g,K)$ module,
$\alpha$ a $K$ type. We say that $\alpha$ is Li positive in $M$ if
\[
   M(\alpha)\neq 0,
\]
and there is a $K$ invariant Hilbert space structure on $M(\alpha)$
so that the matrix coefficient
\[
  c_{M,\alpha} : G\rightarrow \End_\C(M(\alpha))
\]
is Li positive.

Notice that the concept of Li positivity can be defined for
disconnected real reductive groups.

\subsection{Finite dimensional representations and spherical representations}

The following two examples are not related to our Eisenstein
integral. We omit their easy proofs.

\begin{finitedimensional}
Let
\[
   \rho : G\rightarrow \GL(M_{\rho})
\]
be a finite dimensional continuous representation of $G$. Assume
$G$ has a compact center. Then there is a $K$ invariant Hilbert
space structure on $M_\rho$ so that $\rho(x)$ is positive definite
for every $x\in \exp(\p_0)$.
\end{finitedimensional}

\begin{trivial}
If $M$ is a spherical irreducible $(\g,K)$ module with real
infinitesimal character, then the trivial $K$ type is Li positive
in $M$.
\end{trivial}

\subsection{Li positivity of our Eisenstein integrals}

Let $\sigma$ be a finite dimensional unitary representation of $K$
on $W$. Let $\sigma'$ be the representation of $K'$ on $W'$ which
is $\q$-associated to $\sigma$. Notice that $\sigma'$ is also a
unitary representation. The Hilbert space structure on $W'$ is the
restriction of that on $W$. Define a representation $\tau$ of
$K\times K$ on $\End_\C(W)$, a representation $\tau'$ of $K'\times
K'$ on $\End_\C(W')$, and an injection
\[
    j_\q : \End_\C(W')\rightarrow \End_\C(W)
\]
as in section \ref{S:Mlayer}.

\begin{eisensteinli}\label{T:eisensteinli}
Let $\phi\in \con^\infty(G'_\C/U';\tau')$. If $\phi$ is Li
positive, then
\[
   E_\q(j_\q\circ \phi)\in \con^\infty(G_\C/U;\tau)
\]
is also Li positive.
\end{eisensteinli}

\proof
Notice that
\[
   W=W'\oplus\sigma(\bar{\n}\cap \k)W
\]
is an orthogonal decomposition. Therefore $j_\q$ maps positive
definite operators to positive definite operators. Write
\[
  \phi'=j_\q\circ\phi.
\]
Then $\dqr \phi'$ has positive definite values everywhere. Notice
that $\sigma(k^{-1})$ is the adjoint operator of $\sigma(k)$ for all
$k\in K$. Therefore, for every $x\in G_\C/U$,
\begin{eqnarray*}
   &&\quad E_\q(j_\q\circ \phi : x)\\
   &&=\int_K\tau(k^{-1},k^{-1})(\dqr\phi')(H_\q(kx))\,dk\\
   &&=\int_K \sigma(k^{-1})\circ ((\dqr \phi')(H_\q(kx)))\circ \sigma(k)\,dk
\end{eqnarray*}
is positive definite. \qed

\subsection{Li positivity of Bottom layers}

Let $M'$ be a finitely generated admissible $(\g',K')$ module.
$M=\RR^S(M')$ as before. Let $\alpha$ be a $K$ type which is in
the bottom layer of $M$. Let $\alpha'$ be the $K'$ type
$\q$-associated to $\alpha$.

\begin{cohomologicallipositive}\label{T:cohomologicallipositive}
If $\alpha'$ is Li positive in $M'$, then $\alpha$ is Li positive
in $M$.
\end{cohomologicallipositive}

\proof Fix an irreducible unitary representation $\sigma_0$ of $K$
of type $\alpha$ on a vector space $W_\alpha$. Let $\sigma'_0$ be
the representation of $K'$ on $W_{\alpha'}$ which is
$\q$-associated to $\sigma_0$. Then $\sigma'_0$ is an irreducible
unitary representation of $K'$ of type $\alpha'$. Denote by
$\la\,,\,\ra_{\alpha}$ the inner product on $W_\alpha$. The inner
product on $W_{\alpha'}$ is
\[
  \la\,,\,\ra_{\alpha'}=\la\,,\,\ra_{\alpha}|_{W_{\alpha'}\times
  W_{\alpha'}}.
\]

Write
\[
  W'=M'(\alpha')=L\otimes W_{\alpha'}
\]
as before. And write
\[
  W=M(\alpha)=L\otimes W_{\alpha}.
\]
Define $\tau$ and $\tau'$ as before. There is a $K'$ invariant
Hilbert space structure $\la\,,\,\ra'$ on $W'$ so that
\[
   c_{M',\alpha'}\in \con^\infty(G';\tau')
\]
is Li positive. There is a unique Hilbert space structure
$\la\,,\,\ra_L$ on $L$ so that
\[
  \la\,,\,\ra'=\la\,,\,\ra_L \otimes \la\,,\,\ra_{\alpha'}.
\]
Define a $K$ invariant Hilbert structure
\[
    \la\,,\,\ra=\la\,,\,\ra_L \otimes \la\,,\,\ra_{\alpha}
\]
on $W$. Then
\[
  \la\,,\,\ra'=\la\,,\,\ra|_{W'\times W'}.
\]
We conclude the proof by using Lemma \ref{T:liflensted}, Theorem
\ref{T:coef4} and Proposition \ref{T:eisensteinli}. \qed

\subsection{$A_\q(\lambda)$}

We use the notations in section \ref{S:aqlambda}. Assume that
\[
  \Lambda=\lambda|_\t+\delta_{\n\cap\p}|_\t
\]
is dominant with respect to $\n\cap \k$. Let $\a_0$ be the
centralizer of $\t$ in $\p_0$. Notice that $(\t\cap \k_0)\oplus
\a_0$ is a fundamental Cartan subalgebra of $\g_0$. The above
theorem implies that if $\lambda|_{\a_0}$ is real, then the $K$
type $\alpha$ in the bottom layer is Li positive in
$A_\q(\lambda)$. In particular, if $M$ is a discrete series, a
limit of discrete series, or an irreducible unitary representation
with nonzero cohomology, then the unique minimal $K$ type of $M$
is Li positive. The first two cases are essentially known to
Flensted-Jensen \cite{FJ}, and some of the third case are known to
Michael Harris and Jian-Shu Li \cite{HL}.

\subsection{Lowest weight module}

Assume $G$ is of Hermitian type and
\[
  \p=\p^+ \oplus \p^-,
\]
where both $\p^+$ and $\p^-$ are $K$ invariant abelian subspaces of
$\p$, and
\[
  \p^+ = \overline{\p^-}.
\]
Assume
\[
  \q=\k\oplus \p^+.
\]
Fix an irreducible continuous unitary representation $\sigma_0$ of
$K$ on $W_\alpha$, of type $\alpha$. Define a $(\g,K)$ module
\[
  M=\U(\g)\otimes_{\U(\bar \q)}W_\alpha.
\]
Here we view $W_\alpha$ a $\bar \q$ module via the trivial $\p^-$
action. The $\U(\g)$ action on $M$ is via left multiplication. The
$K$ action on $M$ is via the tensor product. It's clear that
\[
  M(\alpha)=W_\alpha.
\]
It' easy to determine the action of $\U(\g)^K$ on $M(\alpha)$.
Then we are able to express $d_{M,\alpha}$ in terms of Eisenstein
integral for $\theta$ stable parabolic. It turns out that
\[
  c_{M,\alpha}(kx^2)=\int_K \sigma_0(kl^{-1} \tilde H_\q
  (lx)l)\,dl, \quad k\in K, \, x\in \exp(\p_0),
\]
where
\[
   \tilde H_\q : G \rightarrow K_\C=G'_\C
\]
is defined in section \ref{S:aqlambda}. Notice that
\[
  \tilde H_\q(G)\subset \exp(\sqrt{-1}\k_0),
\]
and therefore $\sigma_0\circ \tilde H_\q$ has positive definite
values everywhere. This implies that $\alpha$ is Li positive in
$M$. Nollan Wallach and Jian-Shu Li have another expression of
$\Tr\circ c_{\M,\alpha}$ in the unitary case \cite{Li1}.

\subsection{$\SL_2(\R)$}

\begin{sl2r}
Assume $G=\SL_2(\R)$. Let $M$ be an irreducible $(\g,K)$ module.
If $M$ has a real infinitesimal character, then every minimal $K$
type is Li positive in $M$.
\end{sl2r}

\proof We use the well know result on the classification of
irreducible representations of $\SL_2(\R)$. There are 4 cases:
\begin{enumerate}
    \item
      $M$ is a finite dimensional $(\g,K)$ module;
    \item
      $M$ is the $(\g,K)$ module of a discrete series
      representation;
    \item
      $M$ is the $(\g,K)$ module of a limit of discrete series
      representation;
    \item
      $M$ is the $(\g,K)$ module of an irreducible principle series representation.
\end{enumerate}
The first 3 cases are already known. The fourth case can be proved
by an explicit calculation and a trivial estimate. \qed

\subsection{Branching rules}

Let's mention a useful branching rule relative to Li positivity.

Let $M$ be a finitely generated admissible $(\g,K)$ module. Recall
that there is a canonical globalization (or called smooth
globalization) $\widehat M$ of $M$. $\widehat M$ is uniquely
determined by the following
\begin{enumerate}
  \item
    $\widehat M$ is a smooth Frechet representation of G;
  \item
    The $(\g,K)$ module of $\widehat M$ is $M$;
  \item
    $\widehat M$ has moderate growth.
\end{enumerate}
For more details of the globalization, see chapter 11 of
\cite{W2}. We note that every element of $M^{\oc}$ extends to a
continuous linear functional on $\widehat M$. Let $u\in M^{\oc}$,
$v\in \widehat M$. Recall that the matrix coefficient
$c_M(u\otimes v)\in \con^\infty(G;\C)$ is defined by
\[
  c_M(u\otimes v)(g)=u(g.v), \qquad g\in G.
\]

\begin{l1}
Let $M$ be a finitely generated admissible $(\g,K)$ module, $H$ a
unimodular closed subgroup of $G$. Fix a Haar measure on $H$. $M$ is
said to be $\oL^1$ on $H$ if
\[
  c_M(u\otimes v)|_H\in \oL^1(H) \quad \textrm{for every} \quad u\in M^{\oc},\,v\in \widehat M,
\]
and the map
\[
  \begin{array}{rcl}
    \widehat M & \rightarrow & \oL^1(H),\\
             v &\mapsto  &c_M(u\otimes v)|_H
  \end{array}
\]
is continuous for every $u\in M^{\oc}$.
\end{l1}

Let $G_1$ be a closed real reductive subgroup of $G$ so that
$\theta(G_1)=G_1$. Then $\theta$ induces an Cartan involution on
$G_1$. Denote by
\[
  K_1=K\cap G_1,
\]
and $\g_1$ the complexified Lie algebra of $G_1$.

\begin{branching}
Let $M$ be a finitely generated admissible $(\g,K)$ module, $M_1$ a
finitely generated admissible $(\g_1,K_1)$ module. Let $\alpha$ be a
$K$ type and $\alpha_1$ a $K_1$ type. Assume that
\begin{enumerate}
  \item
    The finitely generated $(\g\times \g_1, K\times K_1)$ module
    $M\otimes M_1^{\oc}$ is $\oL^1$ on the subgroup $\Delta(G_1)$ of
    $G\times G_1$, where $\Delta(G_1)$ means the diagonal subgroup;
 \item
    $\alpha$ is Li positive in $M$; $\alpha_1$ is Li positive in $M_1$.
\end{enumerate}
Then
\[
  \dim(\Hom_{(g_1,K_1)}(M,M_1))\geq
  \dim(\Hom_{K_1}(M(\alpha),M_1(\alpha_1))).
\]
\end{branching}

There is also a unitary version of the above branching rule.

\subsection{Final remarks}
Let $M$ be an irreducible $(\g,K)$ module, $\alpha$ a $K$ type
occurring in $M$. Assume now $\alpha$ is a minimal $K$ type. We use
Vogan-Zuckerman's classification of irreducible $( \g,K)$ modules.
There is a $\theta$ stable parabolic subalgebra $\q$, two $(\g',K')$
module $M'$ and $P'$, depending on $M$, such that
\begin{enumerate}
   \item
     $M=\RR^S(M')$, where $S=\dim(\n\cap\k)$ as usual;
   \item
     $G'$ is quasi-split; $P'$ is a principle
     series $(\g',K')$ module;
   \item
     $\alpha$ is in the bottom layer of $M$, $\alpha'$ is a fine
     $K'$ type occurring in $P'$, where $\alpha'$ is the $K'$ type
     $\q$-associated to $\alpha$;
   \item
     $M'$ is the irreducible subquotient of $P'$ so that $M'(\alpha')=P'(\alpha')$.
\end{enumerate}
Theorem \ref{T:coef4} enable us to express $c_{M,\alpha}$ in terms
of $c_{M',\alpha'}$. We know that
\[
   c_{M',\alpha'}=c_{P',\alpha'}
\]
can be easily expressed as an ordinary Eisenstein integral. Theorem
\ref{T:cohomologicallipositive} reduces the question of Li
positivity of $\alpha$ to that of $\alpha'$. But in general, it's
difficult to see when does an ordinary Eisenstein integral preserves
Li positivity. This phenomena is opposite to the question on
unitarizability: it's easy to see that ordinary parabolic inductions
preserves unitarity, but it's rather difficult to see when is the
cohomological induction preserves unitarity. I hope the following is
true for a large class of $G$ (Maybe all):

Let $M$ be an irreducible $(\g,K)$ module with real infinitesimal
character $\lambda$. Let $\alpha$ be a $K$ type occurring in it. If
the size of $\alpha$ is smaller than the size of $\lambda$ in some
sense, then $\alpha$ is Li positive in $M$.

The definition of Eisenstein integral for $\theta$ stable
parabolic subalgebras can be generalized to reductive symmetric
spaces, as Flensted-Jensen did for discrete series.

For applications to branching rules and other aspects of
representation theory, we need a careful study of the growth
condition. On the other hand, our results have some consequences
on the growth of representations. We hope to do these in other
papers.

\end{document}